\documentclass[a4paper,12pt]{article}

\usepackage[all]{xy}
\usepackage[T1]{fontenc}
\usepackage[dvips]{graphicx}
\usepackage{amsmath,amssymb,amsfonts,amsthm,latexsym,mathrsfs,textcomp,verbatim, enumerate, soul, color, MnSymbol}
\xyoption{web}

\title{Gelfand spectra and\\ Wallman compactifications}

\author{Olivia Caramello \vspace{3 mm}\\ {\small DPMMS, University of Cambridge,}\\{\small Wilberforce Road, Cambridge CB3 0WB, U.K.}\\{\small O.Caramello@dpmms.cam.ac.uk}\thanks{The author gratefully acknowledges the support of a Research Fellowship from Jesus College, Cambridge (U.K.)}}

\date{April 15, 2012}
\begin{document}

\mathcode`\<="4268  
\mathcode`\>="5269  
\mathcode`\.="313A  
\mathchardef\semicolon="603B 
\mathchardef\gt="313E
\mathchardef\lt="313C

\newcommand{\app}
 {{\sf app}}

\newcommand{\Ass}
 {{\bf Ass}}

\newcommand{\ASS}
 {{\mathbb A}{\sf ss}}

\newcommand{\Bb}
{\mathbb}

\newcommand{\biimp}
 {\!\Leftrightarrow\!}

\newcommand{\bim}
 {\rightarrowtail\kern-1em\twoheadrightarrow}

\newcommand{\bjg}
 {\mathrel{{\dashv}\,{\vdash}}}

\newcommand{\bstp}[3]
 {\mbox{$#1\! : #2 \bim #3$}}

\newcommand{\cat}
 {\!\mbox{\t{\ }}}

\newcommand{\cinf}
 {C^{\infty}}

\newcommand{\cinfrg}
 {\cinf\hy{\bf Rng}}

\newcommand{\cocomma}[2]
 {\mbox{$(#1\!\uparrow\!#2)$}}

\newcommand{\cod}
 {{\rm cod}}

\newcommand{\comma}[2]
 {\mbox{$(#1\!\downarrow\!#2)$}}

\newcommand{\comp}
 {\circ}

\newcommand{\cons}
 {{\sf cons}}

\newcommand{\Cont}
 {{\bf Cont}}

\newcommand{\ContE}
 {{\bf Cont}_{\cal E}}

\newcommand{\ContS}
 {{\bf Cont}_{\cal S}}

\newcommand{\cover}
 {-\!\!\triangleright\,}

\newcommand{\cstp}[3]
 {\mbox{$#1\! : #2 \cover #3$}}

\newcommand{\Dec}
 {{\rm Dec}}

\newcommand{\DEC}
 {{\mathbb D}{\sf ec}}

\newcommand{\den}[1]
 {[\![#1]\!]}

\newcommand{\Desc}
 {{\bf Desc}}

\newcommand{\dom}
 {{\rm dom}}

\newcommand{\Eff}
 {{\bf Eff}}

\newcommand{\EFF}
 {{\mathbb E}{\sf ff}}

\newcommand{\empstg}
 {[\,]}

\newcommand{\epi}
 {\twoheadrightarrow}

\newcommand{\estp}[3]
 {\mbox{$#1 \! : #2 \epi #3$}}

\newcommand{\ev}
 {{\rm ev}}

\newcommand{\Ext}
 {{\rm Ext}}

\newcommand{\fr}
 {\sf}

\newcommand{\fst}
 {{\sf fst}}

\newcommand{\fun}[2]
 {\mbox{$[#1\!\to\!#2]$}}

\newcommand{\funs}[2]
 {[#1\!\to\!#2]}

\newcommand{\Gl}
 {{\bf Gl}}

\newcommand{\hash}
 {\,\#\,}

\newcommand{\hy}
 {\mbox{-}}

\newcommand{\im}
 {{\rm im}}

\newcommand{\imp}
 {\!\Rightarrow\!}

\newcommand{\Ind}[1]
 {{\rm Ind}\hy #1}

\newcommand{\iten}[1]
{\item[{\rm (#1)}]}

\newcommand{\iter}
 {{\sf iter}}

\newcommand{\Kalg}
 {K\hy{\bf Alg}}

\newcommand{\llim}
 {{\mbox{$\lower.95ex\hbox{{\rm lim}}$}\atop{\scriptstyle
{\leftarrow}}}{}}

\newcommand{\llimr}
 {{\mbox{$\lower.95ex\hbox{{\rm lim}}$}\atop{\scriptstyle
{\rightarrow}}}{}}

\newcommand{\llimd}
 {\lower0.37ex\hbox{$\pile{\lim \\ {\scriptstyle
\leftarrow}}$}{}}

\newcommand{\Mf}
 {{\bf Mf}}

\newcommand{\Mod}
 {{\bf Mod}}

\newcommand{\MOD}
{{\mathbb M}{\sf od}}

\newcommand{\mono}
 {\rightarrowtail}

\newcommand{\mor}
 {{\rm mor}}

\newcommand{\mstp}[3]
 {\mbox{$#1\! : #2 \mono #3$}}

\newcommand{\Mu}
 {{\rm M}}

\newcommand{\name}[1]
 {\mbox{$\ulcorner #1 \urcorner$}}

\newcommand{\names}[1]
 {\mbox{$\ulcorner$} #1 \mbox{$\urcorner$}}

\newcommand{\nml}
 {\triangleleft}

\newcommand{\ob}
 {{\rm ob}}

\newcommand{\op}
 {^{\rm op}}
 
\newcommand{\palrr}[4]{ 
  \def\labelstyle{\scriptstyle} 
  \xymatrix{ {#1} \ar@<0.5ex>[r]^{#2} \ar@<-0.5ex>[r]_{#3} & {#4} } } 
  
\newcommand{\palrl}[4]{ 
  \def\labelstyle{\scriptstyle} 
  \xymatrix{ {#1} \ar@<0.5ex>[r]^{#2}  &  \ar@<0.5ex>[l]^{#3} {#4} } }  

\newcommand{\pepi}
 {\rightharpoondown\kern-0.9em\rightharpoondown}

\newcommand{\pmap}
 {\rightharpoondown}

\newcommand{\Pos}
 {{\bf Pos}}

\newcommand{\prarr}
 {\rightrightarrows}

\newcommand{\princfil}[1]
 {\mbox{$\uparrow\!(#1)$}}

\newcommand{\princid}[1]
 {\mbox{$\downarrow\!(#1)$}}

\newcommand{\prstp}[3]
 {\mbox{$#1\! : #2 \prarr #3$}}

\newcommand{\pstp}[3]
 {\mbox{$#1\! : #2 \pmap #3$}}

\newcommand{\relarr}
 {\looparrowright}

\newcommand{\rlim}
 {{\mbox{$\lower.95ex\hbox{{\rm lim}}$}\atop{\scriptstyle
{\rightarrow}}}{}}

\newcommand{\rlimd}
 {\lower0.37ex\hbox{$\pile{\lim \\ {\scriptstyle
\rightarrow}}$}{}}

\newcommand{\rstp}[3]
 {\mbox{$#1\! : #2 \relarr #3$}}

\newcommand{\scn}
 {{\bf scn}}

\newcommand{\scnS}
 {{\bf scn}_{\cal S}}

\newcommand{\semid}
 {\rtimes}

\newcommand{\Sep}
 {{\bf Sep}}

\newcommand{\sep}
 {{\bf sep}}

\newcommand{\Set}
 {{\bf Set}}

\newcommand{\Sh}
 {{\bf Sh}}

\newcommand{\ShE}
 {{\bf Sh}_{\cal E}}

\newcommand{\ShS}
 {{\bf Sh}_{\cal S}}

\newcommand{\sh}
 {{\bf sh}}

\newcommand{\Simp}
 {{\bf \Delta}}

\newcommand{\snd}
 {{\sf snd}}

\newcommand{\stg}[1]
 {\vec{#1}}

\newcommand{\stp}[3]
 {\mbox{$#1\! : #2 \to #3$}}

\newcommand{\Sub}
 {{\rm Sub}}

\newcommand{\SUB}
 {{\mathbb S}{\sf ub}}

\newcommand{\tbel}
 {\prec\!\prec}

\newcommand{\tic}[2]
 {\mbox{$#1\!.\!#2$}}

\newcommand{\tp}
 {\!:}

\newcommand{\tps}
 {:}

\newcommand{\tsub}
 {\pile{\lower0.5ex\hbox{.} \\ -}}

\newcommand{\wavy}
 {\leadsto}

\newcommand{\wavydown}
 {\,{\mbox{\raise.2ex\hbox{\hbox{$\wr$}
\kern-.73em{\lower.5ex\hbox{$\scriptstyle{\vee}$}}}}}\,}

\newcommand{\wbel}
 {\lt\!\lt}

\newcommand{\wstp}[3]
 {\mbox{$#1\!: #2 \wavy #3$}}
 
\newcommand{\fu}[2]
{[#1,#2]}


%
%
%
\def\pushright#1{{
   \parfillskip=0pt            
   \widowpenalty=10000         
   \displaywidowpenalty=10000  
   \finalhyphendemerits=0      
  %
   \leavevmode                 
   \unskip                     
   \nobreak                    
   \hfil                       
   \penalty50                  
   \hskip.2em                  
   \null                       
   \hfill                      
   {#1}                        
  %
   \par}}                      

\def\qed{\pushright{$\square$}\penalty-700 \smallskip}

\newtheorem{theorem}{Theorem}[section]

\newtheorem{proposition}[theorem]{Proposition}

\newtheorem{scholium}[theorem]{Scholium}

\newtheorem{lemma}[theorem]{Lemma}

\newtheorem{corollary}[theorem]{Corollary}

\newtheorem{conjecture}[theorem]{Conjecture}

\newenvironment{proofs}%
 {\begin{trivlist}\item[]{\bf Proof }}%
 {\qed\end{trivlist}}

  \newtheorem{rmk}[theorem]{Remark}
\newenvironment{remark}{\begin{rmk}\em}{\end{rmk}}

  \newtheorem{rmks}[theorem]{Remarks}
\newenvironment{remarks}{\begin{rmks}\em}{\end{rmks}}

  \newtheorem{defn}[theorem]{Definition}
\newenvironment{definition}{\begin{defn}\em}{\end{defn}}

  \newtheorem{eg}[theorem]{Example}
\newenvironment{example}{\begin{eg}\em}{\end{eg}}

  \newtheorem{egs}[theorem]{Examples}
\newenvironment{examples}{\begin{egs}\em}{\end{egs}}


\bgroup           
\let\footnoterule\relax  
\maketitle

\begin{abstract}
We carry out a systematic, topos-theoretically inspired, investigation of Wallman compactifications with a particular emphasis on their relations with Gelfand spectra and Stone-\v{C}ech compactifications. In addition to proving several specific results about Wallman bases and maximal spectra of distributive lattices, we establish a general framework for functorializing the representation of a topological space as the maximal spectrum of a Wallman base for it, which allows to generate different dualities between categories of topological spaces and subcategories of the category of distributive lattices; in particular, this leads to a categorical equivalence between the category of commutative $C^{\ast}$-algebras and a natural category of distributive lattices. We also establish a general theorem concerning the representation of the Stone-\v{C}ech compactification of a locale as a Wallman compactification, which subsumes all the previous results obtained on this problem.     
\end{abstract} 
\egroup 

\newpage

\tableofcontents

\section{Introduction}

This paper consists in a systematic investigation of Wallman compactifications in relation to Gelfand spectra, Stone-\v{C}ech compactifications and more generally to the representation theory of topological spaces. It is shown that the notion of Wallman base can serve in many contexts as a convenient tool for representing topological spaces, to the point of leading to useful dualities between notable categories of topological spaces, such as the category of $T_{1}$ compact spaces or that of compact Hausdorff spaces, and natural categories of distributive lattices.   

Our analysis is inspired by the view of Grothendieck toposes introduced in \cite{OC10}, which regards toposes as `unifying spaces' being able to effectively act as `bridges' between different representations of a given mathematical object (whenever the latter can be formalized as different sites of definition for one topos). In fact, the notion of base for a topological space is particularly amenable to such an approach, since by Grothendieck's Comparison Lemma for any topological space $X$ and base $B$ for it we have an equivalence of toposes
\[
\Sh(X)\simeq \Sh(B, J_{{\cal O}(X)}^{can}|B),
\]   
where $J_{{\cal O}(X)}^{can}|B$ is the Grothendieck topology induced on $B$ by the canonical topology $J_{{\cal O}(X)}^{can}$ on the frame ${\cal O}(X)$ of open sets of $X$. Such representation, which can be interpreted as an `objectification' of the abstract relationship between $X$ and $B$, indicates that when the topology $J_{{\cal O}(X)}^{can}|B$ can be characterized `intrinsically' in terms of the partially ordered structure of $B$ induced by the inclusion $B\subseteq {\cal O}(X)$, the space $X$ admits an `intrinsic' representation in terms of $B$ which, if appropriately functorialized, can lead to a duality between a category of such spaces $X$ and a category of such posets $B$. For instance, if $X$ is compact and $B$ is a normal Wallman base for it then, under a form of the axiom of choice, the topology $J_{{\cal O}(X)}^{can}|B$ can be identified with a Grothendieck topology $J_{m}^{B}$ intrinsically defined in terms of the lattice structure on $B$, and in the particular case when $X$ is Hausdorff and $B$ is equal to the lattice $Coz(X)$ of co-zero sets on $X$ (as defined in \cite{stone}) the topology $J_{m}^{B}$ admits even a further representation, as the Grothendieck topology on $B$ whose covering sieves are those which contain countable covering families. In fact, this latter choice leads to a duality between the category of compact Hausdorff spaces and particular category of distributive lattices whose objects are Alexandrov algebras satisfying a natural lattice-theoretic condition and whose arrows are the distributive lattices between them which preserve countable joins; composed with Gelfand duality for commutative $C^{\ast}$-algebras, this duality gives rise to a categorical equivalence between the category of $C^{\ast}$-algebras and such category of lattices, which can be exploited for studying $C^{\ast}$-algebras from a purely lattice-theoretic viewpoint. Another situation in which the topology $J_{{\cal O}(X)}^{can}|B$ can be characterized intrinsically is when $X$ is a $T_{1}$ compact space and $B$ is equal to ${\cal O}(X)$; this leads to a duality between the category of $T_{1}$ compact spaces and a particular category of frames, which restricts to a different duality for compact Hausdorff spaces. Of course, different choices of Wallman bases for compact spaces can lead to different dualities, and in fact we present a general framework for generating such dualities in section \ref{generalframework}.

The above-mentioned topos-theoretic viewpoint also guides us in our investigation of the relationships between the Stone-\v{C}ech compactification of a topological space $X$ and its Wallman compactifications. The study of such relations had been initiated by Wallman himself, who proved that for any normal completely regular space $X$ its Stone-\v{C}ech compactification can be identified with the Wallman compactification $Max({\cal O}(X))$, and was continued by several authors, including Gillman and Jerison \cite{GJ}, Frink \cite{frink} and Johnstone \cite{stone} and \cite{wallman}; in section \ref{WSC} we prove a general theorem, based on the concept of $A$-conjunctive sublattice of a frame $A$, which subsumes all the previous results obtained on this problem and allows to establish (iso)morphisms between the Stone-\v{C}ech compactification of a locale and its Wallman compactifications in new many cases which were not covered by the past treatments.    

The careful reader will appreciate the unifying power that the topos-theoretic viewpoint can offer on these questions. In fact, the naturality of this approach is also witnessed by the possibility of naturally interpreting the different constructions of Gelfand spectra and Stone-\v{C}ech compactifications as Morita-equivalences between different geometric (propositional) theories having the same classifying topos, which in fact we analyze and exploit systematically in the course of the paper. Actually, the topos-theoretic viewpoint is the only one which allows to understand and investigate the relationships between different \emph{ways} of constructing a certain topological space in a unified way, that is as different \emph{representations} (namely, sites of definitions) of a single object (namely, the topos associated to the space); the given topos can then be effectively used as a `bridge' for transferring information between its different representations according to the methodologies introduced in \cite{OC10} (in fact, essentially all of the results obtained in the paper arise as applications of this general technique).    

Let us now proceed to describing the contents of the paper in greater detail.

In section \ref{WallmanStone} we make a systematic study of the concept of Wallman base by identifying its natural lattice-theoretic counterpart, namely the concept of $A$-conjunctive sublattice of a frame $A$, and investigate maximal spectra of distributive lattices from both a point-based and point-free perspective; the integration of these different approaches, combined with some central results in Topos Theory, leads us to several concrete results about Wallman bases and conjunctive lattices, notably including the above-mentioned general theorem relating the Stone-\v{C}ech compactification of a locale and its Wallman compactifications. 

In section \ref{dualtopoldlat} we address the problem of functorializing the representations of topological spaces as maximal spectra of Wallman bases for them, and establish a general duality theorem between an appropriate category of topological spaces each of which equipped with a Wallman base on it and a subcategory of the category of distributive lattices. This duality theorem is then applied in section \ref{T1} to generate a duality for $T_{1}$ compact spaces and in section \ref{dualityalexandrovalg} to obtain a duality between the category of compact Hausdorff spaces and a particular category of Alexandrov algebras. This latter duality is also analyzed, in view of Gelfand duality between commutative $C^{\ast}$-algebras and compact Hausdorff spaces, from the point of view of $C^{\ast}$-algebras leading to an explicit categorical equivalence between the category of $C^{\ast}$-algebras and this category of lattices; in particular, any $C^{\ast}$-algebra is shown to be recoverable from the associated Alexandrov algebra through a construction of essentially order-theoretic and arithmetic nature. The results are presented for \emph{real} $C^{\ast}$-algebras (that is, for rings of real-valued continuous functions on a compact Hausdorff space) but they can be straightforwardly extended to the context of complex $C^{\ast}$-algebras.

In section \ref{Gelfandspectra} we investigate the notion of maximal spectrum of a commutative ring with unit from the point of view of the distributive lattice consisting of the compact open sets of its Zariski spectrum. This leads to a logical characterization of the topos of sheaves on such spectrum as the classifying topos of a certain propositional geometric theory which, if the spectrum is sober, axiomatizes precisely the maximal ideals of the ring. Next, we explicitly characterize the class of rings with the property that the corresponding distributive lattice is conjunctive, and remark that any finite-dimensional $C^{\ast}$-algebra enjoys this property; this leads in particular to an explicit algebraic characterization of the lattice of co-zero sets on its spectrum as a distributive lattice presented by generators and relations.    

\vspace{0.4cm}

\textbf{Prerequisites and notation.} Even though the results of this paper are concrete and should be understandable by anyone with a basic knowledge of order theory and topology, a familiarity with Locale Theory and Topos Theory is definitely needed to follow the proofs and appreciate the underlying methodologies. The reader is referred to \cite{stone} (resp. to \cite{MM}) for the background of Locale Theory (resp. of Topos Theory) necessary for understanding the paper. If not indicated otherwise, our terminology is standard and borrowed from \cite{stone}. All the distributive lattices considered in the paper are bounded. We refer the reader to \cite{OC11} for the notion of $J$-ideal (resp. of frame $Id_{J}({\cal C})$ of $J$-ideals, of $J$-prime filter) on a preorder $\cal C$ for a (Grothendieck) coverage $J$ on $\cal C$.

\section{Wallman and Stone-\v{C}ech compactifications}\label{WallmanStone}

\subsection{Wallman bases}

We recall from \cite{stone} (IV2.4) that, given a topological space $X$, a \emph{Wallman base} $B$ for $X$ is a sublattice of of the frame ${\cal O}(X)$ of open sets of $X$ which is a base for the topology and satisfies the property that for any $U\in B$ and $x\in U$ there exists $V\in B$ such that $U\cup V=X$ and $x\notin V$. 

Given a distributive lattice $B$, we denote by $Max(B)$ the set of its maximal ideals. Recall that an ideal of a distributive lattice $B$ is a subset $I\subseteq B$ which is a lower-set and satisfies the property that $0_{B}\in I$ and for any $a,b\in I$, $a\vee b\in I$; an ideal $I$ of $B$ is said to be prime if $I\neq B$ and for any $a,b\in B$, $a\wedge b\in I$ implies that either $a\in I$ or $b\in I$. An ideal $I$ of $B$ is said to be maximal if $I\neq B$ and $I$ is not strictly contained in any ideal of $B$, equivalently if for any $a\notin I$, the ideal on $B$ generated by $I\cup \{a\}$ is the whole $B$. Since every maximal ideal is prime (cf. Corollary I 2.4), we can equip $Max(B)$ with the subspace topology induced by the Zariski topology on the Stone spectrum $Spec(B)$ of $B$; in other words, a base for the topology on $Max(B)$ is given by the subsets of the form ${\cal G}_{b}=\{I\in Max(B) \textrm{ | } b\notin I\}$ for $b\in B$ (cf. section \ref{maxspectrum} below for more details).

Assuming a form of the axiom of choice (the maximal ideal theorem for distributive lattices), one can prove that $Max(B)$ is compact (cf. Lemma II 3.5 \cite{stone}). 

For any topological space $X$ and any sublattice $B$ of ${\cal O}(X)$, we have a continuous map 
\[
\eta^{X}_{B}:X\to Spec(B)
\]
sending any point $x\in X$ to the prime ideal $\{b\in B \textrm{ | } x\notin b\}$ of $B$; when the space $X$ can be unambiguously inferred from the context, we simply denote $\eta^{X}_{B}$ by $\eta_{B}$. 

Notice that the map $\eta^{X}_{B}$ is always open and it factors through the subspace inclusion $Max(B)\hookrightarrow Spec(B)$ if and only if $B$ is a Wallman base for $X$.  

If $X$ is a $T_{0}$-space then the map $\eta^{X}_{B}$ is injective and therefore it is an homeomorphism if and only if it is surjective.

\begin{lemma}\label{dense}
Let $X$ be a topological space and $X_{s}$ its sobrification. Then the image of the universal map $\eta:X\to X_{s}$ is dense in $X_{s}$.
\end{lemma}

\begin{proofs}
We can realize $X_{s}$ and $\eta:X\to X_{s}$, up to isomorphism, as follows. $X_{s}$ is the set of completely prime filters on the frame ${\cal O}(X)$ of open sets of $X$, and the map $\eta:X\to X_{s}$ sends any point $x\in X$ to the filter $\{U\in {\cal O}(X) \textrm{ | } x\in U\}$. The topology on $X_{s}$ has as open sets those of the form $F_{u}=\{P\in X_{s} \textrm{ | } u\in P\}$. Now, if $F_{u}$ is non-empty then $u\neq \emptyset$ (by definition of completely prime filter) and hence there is $x\in X$ such that $x\in u$, equivalently $\eta(x)\in F_{u}$. We can thus conclude that every non-empty open set of $X_{s}$ has non-empty intersection with $Im(\eta)$, in other words $Im(\eta)$ is dense in $X_{s}$, as required. 
\end{proofs}

Recall that a distributive lattice $D$ is said to be \emph{normal} if for any $a,b\in D$ such that $a\vee b=1$ there exist $c,d\in D$ such that $c\vee a=b\vee d=1$ and $c\wedge d=0$. In \cite{almax} Johnstone proved that a Wallman base $B$ for a topological space $X$ is (semi-)normal if and only if the space $Max(B)$ is Hausdorff. We thus obtain the following result. 

\begin{corollary}\label{compWallman}
Let $X$ be a $T_{0}$ compact space with a (semi-)normal Wallman base $B$. Then the map $\eta^{X}_{B}:X\to Max(B)$ is an homeomorphism; in particular, $X$ is Hausdorff. 
\end{corollary}

\begin{proofs}
By Lemma \ref{dense}, the image $Im(\eta^{X}_{B})$ in $Max(B)$ of $X$ under the map $\eta^{X}_{B}$ is dense in $Max(B)$, that is its closure coincides with $Max(B)$. Since $X$ is compact and $\eta^{X}_{B}$ is continuous $Im(\eta^{X}_{B})$ is compact as a subspace of $Max(B)$. But, $Max(B)$ being Hausdorff (since $B$ is (semi-)normal), $Im(\eta^{X}_{B})$ is closed in $Max(B)$ and hence coincides with $Max(B)$; in other words, $\eta^{X}_{B}$ is surjective and hence an homeomorphism.   
\end{proofs}

\begin{theorem}\label{thmprelim}
Let $X$ be a topological space with a Wallman base $B$ such that the space $Max(B)$ is sober. Then every prime ideal $I$ of $B$ which is closed under arbitrary unions in ${\cal O}(X)$ is a maximal ideal of $B$.  
\end{theorem}

\begin{proofs}
By our hypotheses we have a map $\eta_{B}:X\to Max(B)$ sending any $x\in X$ to the ideal $\{b\in B \textrm{ | } x\notin b\}$. Now, if $Max(B)$ is sober then, by the universal property of the sobrification of $X$, we have a unique continuous map $\tilde{\eta_{B}}:X_{s}\to Max(B)$ such that $\tilde{\eta_{B}}\circ \eta=\eta_{B}$. It is easy to verify that $\tilde{\eta_{B}}$ sends any filter $P$ in $X_{s}$ to the ideal $B\setminus (P\cap B)$, and that for any $b\in B$, $\tilde{\eta_{B}}^{-1}({\cal G}_{b})=F_{b}$. From the fact that $\tilde{\eta_{B}}\circ \eta=\eta_{B}$ it thus follows, by invoking Lemma \ref{dense}, that the image of $\eta_{B}$ is dense in $Max(B)$; indeed, for any basic open set ${\cal G}_{b}$ of $Max(B)$, if ${\cal G}_{b}$ is non-empty then, by definition of (maximal) ideal on $B$, $b\neq \emptyset$ and hence, since $F_{b}$ has non-empty intersection with $Im(\eta)$, we have $\emptyset \neq \tilde{\eta_{B}}(F_{b} \cap Im(\eta))\subseteq {\cal G}_{b} \cap Im(\eta_{B})$, as required. 

Now, by the results in \cite{OC11}, the sobrification of $X$ can be identified with the space of points of the locale ${\cal O}(X)\simeq Id_{J_{B}}(B)$, where $J_{B}$ is the Grothendieck topology on $B$ induced by the canonical topology on ${\cal O}(X)$, that is with the set $X_{B}$ of $J_{B}$-prime filters on $B$, endowed with the topology whose basic open sets are the subsets of the form ${\cal F}_{b}=\{P\in X_{B} \textrm{ | } b\in P\}$. Concretely, the homeomorphism $X_{s}\to X_{B}$ sends any filter $P$ in $X_{s}$ to the intersection $P\cap B$, and hence the map $\tilde{\eta_{B}}$ corresponds, under this homeomorphism, to the function sending any filter $P'\in X_{B}$ to the complement $B\setminus P'$. We can thus conclude that the complement of any $J_{B}$-prime filter on $B$ is a maximal ideal on $B$. Let us unravel this into more concrete terms. The Grothendieck topology $J_{B}$ on $B$ has as covering sieves on any object $b\in B$ those sieves on $b$ which contain families of subsets of $b$ in $B$ such that their union in ${\cal O}(X)$ is equal to $b$. A $J_{B}$-prime filter on $B$ is thus a subset $F\subseteq B$ such that $F$ is an upper set in $B$, $0_{B}\notin F$, for any $a,b\in B$, $a\in F$ and $b\in F$ implies that $a\cap b\in F$, and for any family of elements of $B$ whose union in ${\cal O}(X)$ belongs to $F$ then at least one of the elements of the family belongs to $F$. The condition that the complement of every $J_{B}$-prime filter on $B$ should be a maximal ideal on $B$ can thus be reformulated as follows: every subset $I$ of $B$ such that $0_{B}\in I$ and $I$ is a prime lower set in $B$ (in the sense that for any $a,b\in B$, $a\cap b\in I$ if and only if either $a\in I$ or $b\in I$) closed under unions in ${\cal O}(X)$ is a maximal ideal of $B$.    
\end{proofs}

\begin{rmk}
\emph{One might wonder if the sufficient condition for maximality of ideals of $B$ given in the statement of the theorem can hold also under less restrictive assumptions. In fact, this condition is precisely equivalent to the fact that $\eta_{B}$ is a continuous map from the sobrification of $X$ to $Max(B)$; and this is equivalent, under the assumption that $Max(B)$ is sober, to the condition that the function $\eta_{B}:X\to Spec(B)$ takes values in $Max(B)$, which, as we saw, is equivalent to the property of $B$ to be a Wallman base and hence does not hold in general for arbitrary $B$.}
\end{rmk}

\begin{theorem}\label{compactmax}
Let $X$ be a topological space and let $B$ be a Wallman base of $X$ such that the map $\eta_{B}$ is surjective on $Max(B)$. Then the sobrification of $X$ is homeomorphic to the sobrification of $Max(B)$, and if $Max(B)$ is sober then the maximal ideals of $B$ are exactly the prime ideals on $B$ which are closed under arbitrary unions in ${\cal O}(X)$. 
\end{theorem} 
  
\begin{proofs}
The subspace inclusion $Max(B)\hookrightarrow Spec(B)$ induces a geometric inclusion 
\[
i:\Sh(Max(B))\to \Sh(Spec(B)).
\]

But $J_{B}$ clearly contains the coherent topology $J_{B}^{coh}$ on $B$ and hence we have a geometric inclusion $j:\Sh(B, J_{B})\hookrightarrow \Sh(B, J_{B}^{coh})$. On the other hand, the surjection $\tilde{\eta_{B}}:X_{B}\to Max(B)$ induces a geometric surjection $s:\Sh(X_{B})\to \Sh(Max(B))$. Now, since the topos $\Sh(B, J_{B})$ has enough points (being equivalent to $\Sh(X)$), it is equivalent to the topos $\Sh(X_{B})$ of sheaves on its space of points, and under this equivalence and the well-known equivalence $\Sh(Spec(B))\simeq \Sh(B, J_{B}^{coh})$, the composite $i\circ s$ corresponds to $j$. Now, the uniqueness (up to isomorphism) of the surjection-inclusion factorization of a geometric morphism ensures that $s$ must be an equivalence, in other words $\tilde{\eta_{B}}$ yields an homeomorphism to the sobrification of $Max(B)$. Hence, in view of the concrete description of the map given in the proof of Theorem \ref{thmprelim}, if $Max(B)$ is sober then the elements of $Max(B)$ are exactly the complements in $B$ of the filters in $X_{B}$, that is the maximal ideals of $B$ are exactly the prime ideals on $B$ which are closed under arbitrary unions in ${\cal O}(X)$.      

Notice that this proof represents an application of the technique `toposes as bridges' of \cite{OC10}.
\end{proofs}

\subsection{The maximal spectrum of a distributive lattice}\label{maxspectrum}

Let $D$ be a distributive lattice. Then the topos $\Sh(D, J^{coh}_{D})$ of sheaves on $D$ with respect to the coherent topology $J^{coh}_{D}$ on it is equivalent to the topos of sheaves on the Stone spectrum $X_{D}$ of $D$, that is to the topological space whose underlying set is the collection of the prime filters on $D$ and whose topology is generated by the following basic open sets: ${\cal F}_{a}:=\{P\in X_{D} \textrm{ | } a\in P\}$, for $a\in D$. The equivalence of toposes
\[
\Sh(D, J^{coh}_{D}) \simeq \Sh(X_{D})
\] 
can be read frame-theoretically as an equivalence 
\[
Id_{J^{coh}_{D}}(D) \cong {\cal O}(X_{D})
\]
between the frame $Id_{J^{coh}_{D}}(D)$ of ideals of $D$ and the frame ${\cal O}(X_{D})$ of open sets of the spectrum $X_{D}$.

We now consider the proper prime filters in $X_{D}$ which are minimal with respect to the subset-inclusion ordering, that is with respect to the specialization ordering on $X_{D}$ induced by the subterminal topology on it. The prime filters on $D$ are precisely the complements in $\mathscr{P}(D)$ of the prime ideals on $D$ (cf. Proposition I.2.2 \cite{stone}), therefore the minimal proper prime filters in $X_{D}$ correspond to the ideals on $D$ which are proper and maximal among the prime ideals on $D$ with respect to the subset-inclusion ordering on $\mathscr{P}(X_{D})$; notice that these latter ideals coincide, if we assume the prime ideal theorem, precisely with the maximal ideals on $D$, since every maximal ideal is prime (cf. Corollary I.2.4 \cite{stone}). Let us assume this condition and denote by $Max(D)$ the set of such ideals. If we endow $Max(D)$ with the topology with basic open sets those of the form ${\cal M}_{a}:=\{I\in Max(D) \textrm{ | } a\notin I\}$ for $a\in D$ then we obtain a topological space which is homeomorphic (under the complement bijection) to a subspace of the spectrum $X_{D}$ and therefore we obtain a subtopos 
\[
i:\Sh(Max(D))\hookrightarrow \Sh(X_{D}).
\]                                            
Recall from \cite{OC11} that the topos $\Sh(X_{D})$ can be regarded as the classifying topos of the \emph{theory of prime filters on $D$}, that is of the propositional theory ${\mathbb T}_{D}$ whose signature has one atomic proposition $F_{a}$ for each element $a\in {\cal C}$, and whose axioms are the following:
\[
(\top \vdash F_{1});
\]
\[
(F_{a} \vdash F_{b})
\] 
for any $a\leq b$ in $D$, 
\[
(F_{a}\wedge F_{b} \vdash F_{a\wedge b})
\]
for any $a, b \in D$, 
\[
(F_{a\vee b} \vdash F_{a} \vee F_{b})
\]
for any $a, b \in D$.

By the duality theorem in \cite{OC6}, the subtopos $i:\Sh(Max(D))\hookrightarrow \Sh(X_{D})$ thus naturally corresponds to a unique quotient ${\mathbb T}_{m}^{D}$ of $\mathbb T$ classified by the topos $\Sh(Max(D))$.  
On the other hand, using the representation $\Sh(X_{D})\simeq Id_{J^{coh}_{D}}(D)$ of the topos $\Sh(X_{D})$, we obtain the existence of a unique Grothen-\\-dieck topology $J_{m}^{D}$ on $D$ containing $J^{coh}_{D}$ such that the topos $\Sh(Max(D))$ can be represented as $\Sh(D, J_{m}^{D})$. 

We shall now provide an explicit description of the topology $J_{m}^{D}$ as well as an axiomatization of the theory ${\mathbb T}_{m}^{D}$. 

To this end, we observe that in general, given any poset site $({\cal C}, J)$ and a morphism of sites $\xi:({\cal C}, J)\to (F, J^{can}_{F})$, where $F$ is a frame and $J^{can}_{F}$ is the canonical topology on $F$, which induces an equivalence of toposes
\[
\Sh({\cal C}, J)\simeq \Sh(F, J^{can}_{F}),
\]
for any surjective frame homomorphism $k:F\to F'$, inducing a subtopos $\Sh(s):\Sh(F')\hookrightarrow \Sh(F)$, the Grothendieck topology $K_{F'}\supseteq J$ on $\cal C$ such that the canonical inclusion $r:\Sh({\cal C}, K_{F'})\hookrightarrow \Sh({\cal C}, J)$ makes the following diagram commute can be described concretely as follows.
\[  
\xymatrix {
\Sh({\cal C}, J) \ar[r]^{\simeq} & \Sh(F)\\
\Sh({\cal C}, K_{F'}) \ar[u]^{r}  \ar[r]_{\simeq} & \Sh(F') \ar[u]^{\Sh(s)}}
\]  
A sieve $S:=\{c_{i}\leq c \textrm{ | } \in I\}$ on $c\in D$ is $K_{F'}$-covering if and only if $k(\mathbin{\mathop{\textrm{\huge $\vee$}}\limits_{i\in I}} \xi(c_{i}))=k(c)$.

Given a distributive lattice $D$, Proposition 1.6 \cite{almax} provides, under the assumption that the maximal ideal theorem holds, an explicit description of the nucleus $k_{D}$ on the frame $Id_{J^{coh}_{D}}(D)$ corresponding to the subtopos $i:\Sh(Max(D))\hookrightarrow \Sh(X_{D})\simeq \Sh(Id_{J^{coh}_{D}}(D))$: for any $I\in Id_{J^{coh}_{D}}(D)$,
\[
k_{D}(I)=\{d\in D \textrm{ | } (\forall b\in D)(d\vee b=1 \imp (\exists c\in I)(b\vee c=1))\}.
\]
In order to apply our characterization of the topology $K_{F'}$ to the particular case of the surjective frame homomorphism $k:Id_{J^{coh}_{D}}(D) \to ({Id_{J^{coh}_{D}}(D)})_{k}$ and the canonical morphism of sites $\xi:(D, J_{D}^{coh})\to (Id_{J^{coh}_{D}}(D), J^{can}_{Id_{J^{coh}_{D}}(D)})$ we observe that in this case we have $\mathbin{\mathop{\textrm{\huge $\vee$}}\limits_{i\in I}}\xi(c_{i})=I_{S}$, where $I_{S}$ is the ideal in $D$ generated by the set $\{c_{i} \textrm{ | } i\in I\}$, equivalently the set 
\[
\{b\in D \textrm{ | } b\leq \mathbin{\mathop{\textrm{\huge $\vee$}}\limits_{i\in J}}c_{i} \textrm{ for some finite set } J\subseteq I\}.
\] 
We thus obtain the following description of the topology $J_{m}^{D}$ (under the assumption that the maximal ideal theorem holds): a sieve $S:=\{c_{i}\leq c \textrm{ | } \in I\}$ on $c\in D$ is $J_{m}^{D}$-covering if and only if for every $d\in D$, $c\vee d=1$ implies that there exists a finite subset $J\subseteq I$ such that $\mathbin{\mathop{\textrm{\huge $\vee$}}\limits_{i\in J}}c_{i} \vee d=1$ (cf. also Exercise II 3.5 \cite{stone}). 

Given a distributive lattice $D$, we call the Grothendieck topology $J_{m}^{D}$ on $D$, also denoted by $J_{m}^{D}$, the \emph{maximal topology} on $D$. We shall say that a distributive lattice homomorphism  $f:D\to D'$ is \emph{maximal} if $Spec(f):Spec(D')\to Spec(D)$ restricts to a continuous map $Max(f):Max(D')\to Max(D)$. We can naturally characterize the maximal homomorphisms in terms of the maximal topologies on the two lattices as follows.

\begin{proposition}\label{maxmorphisms}
Assume the maximal ideal theorem for distributive lattices. Let $f:D\to D'$ be a homomorphism of distributive lattices. If $f$ is a maximal homomorphism then $f$ is a morphism of sites $(D, J_{m}^{D})\to (D', J_{m}^{D'})$ (i.e., for any a sieve $S:=\{c_{i}\leq c \textrm{ | } \in I\}$ on $c\in D$ with the property that for every $d\in D$, $c\vee d=1$ in $D$ implies that there exists a finite subset $J\subseteq I$ such that $\mathbin{\mathop{\textrm{\huge $\vee$}}\limits_{i\in J}}c_{i} \vee d=1$ we have that for every $d'\in D'$, $f(c)\vee d'=1$ in $D'$ implies that there exists a finite subset $K\subseteq I$ such that $\mathbin{\mathop{\textrm{\huge $\vee$}}\limits_{i\in K}}f(c_{i}) \vee d'=1$). If the space $Max(D)$ then also the converse holds.
\end{proposition}

\begin{proofs}
Let $f:D\to D'$ be a morphism of distributive lattices; $f$ induces, as a morphism of sites $(D', J_{D'}^{coh}) \to (D, J_{D}^{coh})$, a geometric morphism 
\[
\Sh(D', J_{D'}^{coh})\to \Sh(D, J_{D}^{coh})
\]
which corresponds, under the canonical equivalences 
\[
\Sh(D, J_{D}^{coh}) \simeq \Sh(Spec(D))
\]
and 
\[
\Sh(D', J_{D'}^{coh}) \simeq \Sh(Spec(D')),
\] 
to the geometric morphism $\Sh(f^{-1}):\Sh(Spec(D'))\to \Sh(Spec(D))$ induced by the continuous map $Spec(f)=f^{-1}:Spec(D')\to Spec(D)$. Therefore if $f$ is maximal then $f^{-1}|:Max(D')\to Max(D)$ induces a geometric morphism 
\[
\Sh(f^{-1}|):\Sh(Max(D'))\to \Sh(Max(D))
\]
such that the morphism $\Sh(f^{-1}|):\Sh(D', J^{D'}_{m})\to \Sh(D, J_{m}^{D})$ corresponding to it under the equivalences
\[
\Sh(D', J^{D'}_{m}) \simeq \Sh(Max(D'))
\]
and 
\[
\Sh(D, J^{D}_{m}) \simeq \Sh(Max(D))
\]
makes the diagram
\[  
\xymatrix {
D \ar[d]^{y'^{D}} \ar[r]^{f} & D' \ar[d]^{y'^{D'}} \\
\Sh(D, J^{D}_{m})  \ar[r]_{{\Sh(f^{-1}|)}^{\ast}} & \Sh(D', J^{D'}_{m}).}
\] 
commute, where 
\[
y'^{D}:D\to \Sh(D, J^{D}_{m})
\]
(resp. $y'^{D'}:D'\to \Sh(D', J^{D'}_{m})$) is the composite of the Yoneda embedding $D\hookrightarrow [D^{\textrm{op}}, \Set]$ (resp. $D'\hookrightarrow [D'^{\textrm{op}}, \Set]$) with the associated sheaf functor $[D^{\textrm{op}}, \Set] \to \Sh(D, J^{D}_{m})$ (resp. $[D'^{\textrm{op}}, \Set] \to \Sh(D', J^{D'}_{m})$). Therefore $f$ is a morphism of sites $(D, J^{D}_{m}) \to (D', J^{D'}_{m})$.

Conversely, it is clear that if $f$ is a morphism of sites $(D, J^{D}_{m}) \to (D', J^{D'}_{m})$ and the space $Max(D)$ is sober then $f^{-1}$ restricts to a map $Max(D')\to Max(D)$.  
\end{proofs}

It is immediate to see, under the hypothesis of the maximal ideal theorem, that a distributive lattice $D$ embeds into ${\cal O}(Max(D))\cong Id_{J_{m}^{D}}(D)$ through the map $\eta_{D}:D\to {\cal O}(Max(D))$ (equivalently, the maximal topology $J_{m}^{D}$ on $D$ is subcanonical) if and only if for any sieve $S:=\{c_{i}\leq c \textrm{ | } \in I\}$ on $c\in D$ with the property for every $d\in D$, $c\vee d=1$ implies that there exists a finite subset $J\subseteq I$ such that $\mathbin{\mathop{\textrm{\huge $\vee$}}\limits_{i\in J}}c_{i} \vee d=1$, $S$ is covering in $D$ in the sense that for any $d\in D$ such that $d\geq c_{i}$ for all $i\in I$ then $d\geq c$. By Proposition 4 \cite{wallman}, one can equivalently characterize these lattices as the conjunctive lattices (in the sense of \cite{wallman} and \cite{HS2}), that is as the distributive lattices $D$ such that for any $a,b\in D$, if for any $c\in D$ such that $a\vee c=1$, $b\vee c=1$ then $a\leq b$.   

It is well-known that a distributive lattice $D$ can always be recovered, up to isomorphism, from the topos $\Sh(D, J_{coh})\simeq \Sh(Spec(D))$, in which it embeds full and faithfully (as its full subcategory of compact subterminals). It is therefore natural to ask whether a similar result holds for conjunctive distributive lattices with respect to the topos $\Sh(D, J_{m}^{D})\simeq \Sh(Max(D))$, that is if it is possible to recover a conjunctive lattice $D$ up to isomorphism from the topos $\Sh(Max(D))$ through some kind of topos-theoretic invariant, at least if $D$ is normal. The answer to this question is negative in general. In fact, any normal Wallman base of a compact Hausdorff space $X$ is a conjunctive distributive lattice $B$ such that $\Sh(X)\simeq \Sh(B, J_{m}^{B})$; but there can be in general different Wallman bases for such a space (for example both the lattice $Coz(X)$ of co-zero sets on $X$ and ${\cal O}(X)$ are). Anyway, any countably compact Alexandrov algebra $D$ can be recovered up to isomorphism from the topos $\Sh(Max(D))$ (as $Coz(Max(D))$) (see Theorem \ref{dualalex} below). In connection with the representation of a maximal spectrum in terms of different distributive lattices, it is worth to remark the following result, which implies in particular that if a surjective homomorphism $f:D\to D'$ between two conjunctive distributive lattices induces an equivalence of toposes $\Sh(f):\Sh(Max(D'))\to \Sh(Max(D))$ then $f$ is an isomorphism.    

\begin{lemma}\label{sur}
Let us assume the maximal ideal theorem. Let $D$ and $D'$ be distributive lattices, of which $D$ conjunctive, and let $f:D\to D'$ be a maximal homomorphism between them. If
\[
\Sh(Max(f)):\Sh(Max(D')) \to \Sh(Max(D'))
\] 
is an equivalence (for example if $Max(f):Max(D')\to Max(D)$ is an homeomorphism) then $f$ is injective.
\end{lemma}

\begin{proofs}
By Proposition \ref{maxmorphisms}, $f:D\to D'$ is a morphism of sites $(D, J_{m}^{D})\to (D', J_{m}^{D'})$ and hence it can be recovered, up to isomorphism, as the restriction of the inverse image $\Sh(Max(f))^{\ast}$ of $\Sh(Max(f))$ to $D\hookrightarrow {\cal O}(Max(D))$. From this representation of $f$ our thesis follows immediately. 
\end{proofs}

Now that we have an explicit description of the Grothendieck topology $J_{m}^{D}$, we can obtain an axiomatization for the theory ${\mathbb T}_{m}^{D}$ introduced above: the quotient ${\mathbb T}_{m}^{D}$ is obtained by adding to the axioms of ${\mathbb T}_{D}$ all the sequents of the form
\[
(F_{c} \vdash \mathbin{\mathop{\textrm{\huge $\vee$}}\limits_{i\in I}}F_{c_{i}})
\]
for any $J$-covering sieve $\{c_{i} \to c \textrm{ | } i\in I\}$ in $D$. 

Notice that if the space $Max(D)$ is sober then the models of the theory ${\mathbb T}_{m}^{D}$ can be identified precisely with the complements in $\mathscr{P}(D)$ of the maximal ideals on $D$; specifically, the following theorem holds.

\begin{theorem}\label{dlatmax}
Assume the maximal ideal theorem for distributive lattices. Let $D$ be a distributive lattice such that the space $Max(D)$ is sober. Then, for any prime ideal $P$ of $D$, $P$ is maximal if and only if for any sieve $S:=\{c_{i}\leq c \textrm{ | } i\in I\}$ on $c\in D$ with the property that for every $d\in D$, $c\vee d=1$ implies that there exists a finite subset $J\subseteq I$ such that $\mathbin{\mathop{\textrm{\huge $\vee$}}\limits_{i\in J}}c_{i} \vee d=1$, $c\notin P$ implies that $c_{i}\notin P$ for some $i\in I$. 
\end{theorem}\qed 

Notice that, since the topos $\Sh(Max(D))$ has enough points, the theory ${\mathbb T}_{m}^{D}$ can alternatively be described as the set of sequents over the signature of ${\mathbb T}_{D}$ which are satisfied by all the complements in $\mathscr{P}(D)$ of the maximal ideals on $D$.

\subsection{Conjunctive lattices}

Given a distributive lattice $D$, let us denote by $\eta_{D}:D\to {\cal O}(Max(D))$ sending to any $d\in D$ the basic open set $\{M\in Max(D) \textrm{ | } d\notin M\}$.  

Conjunctive lattices are a natural lattice-theoretic analogue of Wallman bases; indeed, the following propositions holds. 

\begin{proposition}
Let $D$ a distributive lattice. Then the sublattice $Im(\eta_{D})$ of ${\cal O}(Max(D))$ given by the image of the map $\eta_{D}:D\to {\cal O}(Max(D))$, that is the lattice consisting of the open sets $\{M\in Max(D) \textrm{ | } d\notin M\}$ (for $d\in D$) is a Wallman base of the topological space $Max(D)$, and, under the maximal ideal theorem, $\eta_{D}$ defines an isomorphism $D\cong Im(\eta_{D})$ if and only if $D$ is conjunctive. Conversely, any Wallman base for a topological space is a conjunctive lattice. 
\end{proposition}

\begin{proofs}
Let us prove that $Im(\eta_{D})$ is a Wallman base for $Max(D)$. Given $d\in D$ and $M\in Max(D)$ such that $d\notin M$ we want to show that there exists $e\in D$ such that $e\in M$ and for any $N\in Max(D)$ either $e\notin N$ or $d\notin N$.
Since $M$ is a maximal ideal of $D$ the ideal on $D$ generated by $M\cup \{d\}$ contains $1$, in other words $1=d\vee e$ for some $e\in M$. Clearly, such an element $e$ satisfies the required property, since, as we saw in section \ref{maxspectrum}, $D$ is conjunctive if and only if $\eta_{D}$ is injective, if and only if it yields an isomorphism onto its image. 

Conversely, let $B$ be a Wallman base of a topological space $X$. Let $U, U'\in B$; we want to prove that if for any $Z\in B$, $Z\cup U=X$ implies $Z\cup U'=X$ then $U\subseteq U'$. To show that $U\subseteq U'$ we verify that for any $x\in U$, $x\in U'$. So, let us suppose that $x\in U$. By definition of Wallman base there exists $V\in B$ such that $U\cup V=X$ and $x\notin V$; but then $U'\cup V=X$ whence $x\in U'$.   
\end{proofs}

More generally, given a frame $A$ and a sublattice $B$ of $A$, we define $B$ to be \emph{$A$-conjunctive} if for any $a\in A$ and $b\in B$ if for any $c\in B$, $c\vee b=1$ in $A$ implies $c\vee d=1$ in $A$ for some element $d\in B$ such that $d\leq a$ then $b\leq a$. For any topological space $X$ and Wallman base $B$ for $X$, $B$ is ${\cal O}(X)$-conjunctive; indeed, for any $U\in B$ and $Z\in {\cal O}(X)$, for any $x\in U$ there exists $V\in B$ such that $U\cup V=X$ and $x\notin V$; therefore $Z\cup V=1$ whence $x\in Z$. Notice also that a locale $A$ is subfit (in the sense of \cite{wallman}) if and only if, as a sublattice of itself considered as a frame, it is $A$-conjunctive.

Notice that if $B$ is a sublattice of a frame $A$ which is a base for it then the Comparison Lemma yields an equivalence $\Sh(A, J_{A}^{can})\simeq \Sh(B, J_{A}^{can}|_{B})$; also, since $J^{coh}_{B}$ is contained in $J_{A}^{can}|B$ we have a geometric inclusion 
\[
i_{B}:\Sh(A, J_{A}^{can})\hookrightarrow \Sh(B, J^{coh}_{B})\simeq \Sh(Id(B)).
\]
On the other hand, we have the canonical geometric inclusion 
\[
m_{B}:\Sh(B, J_{m}^{B})\hookrightarrow \Sh(B, J_{B}^{coh}),
\]
which if (and only if) $B$ is normal admits a one-sided inverse 
\[
r_{B}:\Sh(B, J_{B}^{coh})\to \Sh(B, J_{m}^{B})
\]
such that $r_{B}\circ m_{B}\cong 1_{\Sh(B, J_{m}^{B})}$ and the direct image of $r_{B}$ is isomorphic to the associated sheaf functor $\Sh(B, J_{B}^{coh})\to \Sh(B, J_{m}^{B})$ (cf. Proposition 3 \cite{wallman}). 

The following characterization result holds.

\begin{proposition}\label{charconjun}
Let $A$ be a frame and $B$ be a sublattice of $A$ which is a base for $A$. Then the following conditions are equivalent.
\begin{enumerate}[(i)]
\item $B$ is $A$-conjunctive;

\item The map $A\to Id_{J_{m}^{B}}(B)$ sending any element $a\in A$ to the $J_{m}^{B}$-closure of the ideal $\{b\in B \textrm{ | } b\leq a\}$ is injective;

\item The geometric inclusion $i_{B}:\Sh(A, J_{A}^{can})\hookrightarrow \Sh(B, J^{coh}_{B})$ factors (necessarily uniquely up to isomorphism) through the geometric inclusion $m_{B}:\Sh(B, J_{m}^{B})\hookrightarrow \Sh(B, J_{B}^{coh})$. 

\end{enumerate}
Moreover, if $B$ is normal condition $(iii)$ is equivalent to the requirement that the composite geometric morphism $r_{B}\circ i_{B}:\Sh(A, J_{A}^{can}) \to \Sh(B, J_{m}^{B})$ be an inclusion.

\end{proposition}

\begin{proofs}
$(i)\biimp (ii)$. Recall that for any ideal $I$ on $B$, the $J_{m}^{B}$-closure $c_{J_{m}^{B}}(I)$ of $I$ is equal to $\{x\in B \textrm{ | } (\forall b\in B)(x\vee b=1) \imp (\exists y\in I)(y\vee b=1)\}$. For any $a\in A$, let us denote by $I_{a}$ the ideal $\{b\in B \textrm{ | } b\leq a\}$. 
It is immediate to see that for any $x\in B$ and $a\in A$, $x\in c_{J_{m}^{B}}(I_{a})$ if and only if for any $c\in B$, $c\vee x=1$ in $A$ implies $c\vee d=1$ in $A$ for some element $d\in B$ such that $d\leq a$. So the condition of $B$ to be $A$-conjunctive can be reformulated as the requirement that for any $b\in B$ and $a\in A$, $b\in c_{J_{m}^{B}}(I_{a})$ (if and) only if $b\leq a$ in $A$. On the other hand, the map $a\to c_{J_{m}^{B}}(I_{a})$ is injective if and only if for any $a, a'\in A$ $c_{J_{m}^{B}}(I_{a})\subseteq c_{J_{m}^{B}}(I_{a'})$ implies $a\leq a'$ (since for any $a,a'\in A$, $c_{J_{m}^{B}}(I_{a\wedge a'})=c_{J_{m}^{B}}(I_{a})\cap c_{J_{m}^{B}}(I_{a'})$), and this condition is clearly equivalent to demanding that if for any $b\in B$, $b\in c_{J_{m}^{B}}(I_{a})$ implies $b\in c_{J_{m}^{B}}(I_{a'})$ then $a\leq a'$. Therefore, since $B$ is a base for $A$, its $A$-conjunctivity implies the injectivity of the map $a\to c_{J_{m}^{B}}(I_{a})$. Conversely, we want to show that if this map is injective then $B$ is $A$-conjunctive. Suppose that $b\in c_{J_{m}^{B}}(I_{a})$; then $c_{J_{m}^{B}}(I_{b})\subseteq c_{J_{m}^{B}}(I_{a})$ and the injectivity of the map gives $b\leq a$ as required.  

$(i)\biimp (iii)$. The geometric inclusion $i_{B}:\Sh(A, J_{A}^{can})\hookrightarrow \Sh(B, J^{coh}_{B})$ factors (necessarily uniquely up to isomorphism) through the geometric inclusion $m_{B}:\Sh(B, J_{m}^{B})\hookrightarrow \Sh(B, J_{B}^{coh})$ if and only if for every element $a\in A$ the ideal $I_{a}:=\{b\in B \textrm{ | } b\leq a\}$ is $J_{m}^{B}$-closed. But this condition is easily seen to be equivalent to saying that $B$ is $A$-conjunctive. 

Finally, let us prove, under the assumption that $B$ is normal, the equivalence of condition $(iii)$ with the requirement that the composite geometric morphism $r_{B}\circ i_{B}:\Sh(A, J_{A}^{can}) \to \Sh(B, J_{m}^{B})$ be an inclusion. If the composite geometric morphism $r\circ i_{B}:\Sh(A, J_{A}^{can}) \to \Sh(B, J_{m}^{B})$ is an inclusion then the map $a\to c_{J_{m}^{B}}(I_{a})$, which can be identified with the restriction to subterminals of the direct image of this morphism, is injective and the composite of it with the canonical inclusion $Id_{J_{m}^{B}}(B)\subseteq Id(B)$ is equal to the restriction to subterminals $a\to I_{a}$ of the geometric morphism $i_{B}:\Sh(A, J_{A}^{can})\hookrightarrow \Sh(B, J^{coh}_{B})$. Indeed, the $A$-conjunctivity of $B$ implies that for any $a\in A$ and any $b\in B$ if $b\in c_{J_{m}^{B}}(I_{a})$ then $b\leq a$ (equivalently, $b\in I_{a}$), that is for any $a\in A$, $c_{J_{m}^{B}}(I_{a})=I_{a}$. Therefore $m_{B}\circ r\circ i_{B}\cong i_{B}$. Conversely, if $f$ is a geometric inclusion $\Sh(A, J_{A}^{can}) \to \Sh(B, J_{m}^{B})$ such that $m_{B}\circ f\cong i_{B}$ then $f=1\circ f=(r_{B}\circ m_{B})\circ f=r_{B}\circ (m_{B}\circ f)=r_{B}\circ i_{B}$, from which it follows that $r_{B}\circ i_{B}$ is an inclusion (it being equal to $f$), as required.   
\end{proofs}

Thanks to the notion of $A$-conjunctive base we can establish a `point-free version' of Corollary \ref{compWallman} above.

\begin{theorem}\label{compact}
Let $A$ be a frame and $B$ be a sublattice of $A$ which is a base for it. Then the composite geometric morphism $r_{B}\circ i_{B}:\Sh(A, J_{A}^{can}) \to \Sh(B, J_{m}^{B})$ is an equivalence if and only if $B$ is $A$-conjunctive and $A$ is compact. Moreover, if $A$ is compact and $B$ is $A$-conjunctive, $J_{m}^{B}=J_{A}^{can}|B$.  
\end{theorem}

\begin{proofs}
By Proposition \ref{charconjun}, $r_{B}\circ i_{B}$ is an inclusion if and only if $B$ is $A$-conjunctive so it remains to verify that it is an equivalence if and only if $A$ is moreover compact. Clearly, $r_{B}\circ i_{B}$ is an equivalence if and only if the restriction to subterminals of its direct image is a frame isomorphism. We recall from the proof of Proposition \ref{charconjun} that such map is the function $g:A \to Id_{J_{m}^{B}}(B)$ sending any element $a\in A$ to the closure $c_{J_{m}^{B}}(I_{a})$ of the ideal $I_{a}=\{b\in B \textrm{ | } b\leq a\}$. By Proposition \ref{charconjun}, the function $g$ is injective if and only if $B$ is $A$-conjunctive so we have to show that it is surjective if and only if $A$ is compact. Now, clearly $g$ is surjective if and only if for any $J_{m}^{B}$-ideal $I$ on $B$, $I$ contains as element its join $\bigvee I$ in $A$. If $A$ is compact then $\bigvee I \in c_{J_{m}^{B}}(I)=I$ since for every $b\in B$, $b\vee \bigvee I=1$ implies $b\vee d=1$ for some $d\in I$. Conversely, if $g$ is surjective then for any ideal $I$ of $B$, $\bigvee I\in c_{J_{m}^{B}}(I)$, from which it follows that if $\bigvee I=1$ then $1\in I$; that is, $B$ being a base for $A$, $A$ is compact.   

The last part of the theorem follows at once from the fact that since $B$ is $A$-conjunctive then the topology $J_{m}^{B}$ is subcanonical and hence equal to the Grothendieck topology induced on $B$ by the canonical topology on the topos $\Sh(B, J_{m}^{B})$, since we have an equivalence $r_{B}\circ i_{B}:\Sh(A, J_{A}^{can}) \to \Sh(B, J_{m}^{B})$.   
\end{proofs}

\subsection{Wallman vs. Stone-\v{C}ech}\label{WSC}

In this section our aim is to investigate the relationship between the Stone-\v{C}ech compactification $\beta(A)$ of a frame $A$ and the `Wallman compactification' $Id_{J_{m}^{B}}(B)$ for a sublattice $B$ of $A$ which is a base for it.  

Let us recall from \cite{wallman} (Lemma 1(ii) and Proposition 2) that for any distributive lattice $B$ the locale $Id_{J_{m}^{B}}(B)$ is compact and it is regular if and only if $B$ is semi-normal (i.e., for any $a,b\in B$ such that $a\vee b=1$ there exist $c,d\in B$ such that $c\vee a=b\vee d=1$ and $c\wedge d$ belongs to the $J_{m}^{B}$-closure $\overline{(0)}^{J_{m}^{B}}$ of the principal ideal $(0)$). Notice that for any conjunctive distributive lattice $B$, $B$ is normal if and only if it is semi-normal.

Given a frame $A$ and a sublattice $B$ of $A$ which is a base for $A$, we saw above that we have a locale map $g:A\to Id_{J_{m}^{B}}(B)$ whose direct image is the function sending any element $a\in A$ to the $J_{m}^{B}$-closure of the ideal $I_{a}=\{b\in B \textrm{ | } b\leq a\}$. If $B$ is normal, let us denote by $r:Id(B)\to Id_{J_{m}^{B}}(B)$ the retract locale map of Proposition 3 \cite{wallman} and by $i:A\to Id(B)$ the locale map whose direct image is the function sending any $a\in A$ to the ideal $I_{a}$.

Let us thus suppose that $A$ is a frame with a normal sublattice $B$ which is a base for it. Then the locale $Id_{J_{m}^{B}}(B)$ is compact and regular and hence, assuming the axiom of dependent choices, by the universal property of the Stone-\v{C}ech compactification, the locale map $g:A\to Id_{J_{m}^{B}}(B)$ factors uniquely through the canonical map $\eta_{A}:A\to \beta(A)$ as a locale map $h:\beta(A) \to Id_{J_{m}^{B}}(B)$ such that $g=h\circ \eta_{A}$. We can represent this in the following commutative diagram:

\[  
\xymatrix {
A \ar[d]^{\eta_{A}} \ar[dr]^{g} \ar[r]^{i} & Id(B) \ar[d]^{r} \\
\beta(A) \ar[r]^{h} & Id_{J_{m}^{B}}(B).}
\] 

When is the map $h$ an isomorphism? If $B$ is conjunctive we can give a particularly natural answer to this question. Indeed, if $B$ is conjunctive then $B$ embeds canonically into $Id_{J_{m}^{B}}(B)$ and the image of $B$ inside $Id_{J_{m}^{B}}(B)$, that is the subset of principal ideals on $B$, can be identified with the image of $B\subseteq A$ in $Id_{J_{m}^{B}}(B)$ under the direct image of the map $g$; therefore if $h$ is an isomorphism then, by the commutativity of the diagram above, the direct image of $\eta_{A}$ sends $B$ injectively to $\beta(A)$, and the image $\eta_{A}(B)$ is a sublattice of $\beta(A)$ which forms a $\beta(A)$-conjunctive base of it; moreover, by Theorem \ref{compact}, the induced topology $J_{\beta(A)}^{can}|\eta_{A}(B)$ is equal to $J_{m}^{\eta_{A}(B)}$. Conversely, if the direct image of $\eta_{A}$ sends $B$ injectively to a $\beta(A)$-conjunctive base for $\beta(A)$ then $\beta(A)\cong Id_{J_{\beta(A)}^{can}|\eta_{A}(B)}(\eta_{A}(B))$ and $h$ is an isomorphism. 

Notice that from the commutativity of the above diagram it follows that if $g$ is an embedding (equivalently, by Proposition \ref{charconjun}) then $\eta_{A}$ is an embedding (equivalently, by Theorem IV2.1 \cite{stone}, $A$ is completely regular). We record this remark in the following proposition.

\begin{proposition}
Let $A$ be a locale with a normal base $B$ which is $A$-conjunctive. Then $A$ is completely regular. 
\end{proposition}\qed   

Notice also that if $A$ is completely regular then the direct image of $\eta_{A}$ is injective and therefore it sends $B$ bijectively to a subset of $\beta(A)$.

Summarizing, we have the following result.

\begin{theorem}\label{genwallman}
Assume the axiom of dependent choices. Let $A$ be a locale, and $\eta_{A}:A\to \beta(A)$ its Stone-\v{C}ech compactification. Let $B$ be a normal conjunctive sublattice of $A$. Then the map $h:\beta(A)\to Id_{J_{m}^{B}}(B)$ defined above is an isomorphism if and only if the direct image of $\eta_{A}$ sends $B$ injectively to $\beta(A)$ and $\eta_{A}(B)\subseteq \beta(A)$ is a $\beta(A)$-conjunctive base for $\beta(A)$. 

In particular, if $A$ is a completely regular locale and a $B$ is a conjunctive normal base for it (for example, if $B$ is a normal $A$-conjunctive base for $A$) then $h$ is an isomorphism if and only if $\eta_{A}(B)\subseteq \beta(A)$ is a $\beta(A)$-conjunctive base for $\beta(A)$.    
\end{theorem}\qed                                                                                            

\begin{remark}\label{finitejoins}
In the presence of the axiom of choice, for any completely regular locale $L$ the direct image of the universal locale map $\eta_{L}:L\to \beta(L)$ preserves finite joins. Indeed, $\beta(L)$ can be identified with the space ${\cal M}_{L}$ of maximal completely regular filters on $L$, endowed with the topology whose open sets are the subsets of the form ${\cal M}_{L}^{u}:=\{M\in {\cal M}_{L} \textrm{ | } u\in M\}$ for $u\in L$ (cf. \cite{stone}, \cite{BanC} and \cite{BMSC3}), while ${\eta_{L}}_{\ast}:L\to \beta(L)$ can be identified with the map sending any $u\in L$ to ${\cal M}_{L}^{u}$; and ${\eta_{L}}_{\ast}$ preserves finite joins since every maximal regular filter on $L$ is prime (cf. the proof of Proposition IV 2.3 \cite{stone}). 
\end{remark}
           
Theorem \ref{genwallman} represents a significant generalization of Wallman's result (cf. \cite{Wallman}) that for any normal completely regular space $X$ the Stone-\v{C}ech compactification $\beta(X)$ is homeomorphic to the Wallman compactification $Max({\cal O}(X))$. Our result subsumes both Theorem IV 2.7 (originally due to Frink \cite{frink}) and Proposition 8 \cite{wallman}. Indeed, from the theorem we can immediately deduce the fact that the Stone-\v{C}ech compactification of a subfit normal locale $A$ is isomorphic to $Id_{J_{m}^{A}}(A)$ (Proposition 8 \cite{wallman}), as well as the following topological result. 

Below, given a continuous map $f:X\to Y$ of topological spaces, we denote by $f_{\ast}:{\cal O}(X)\to {\cal O}(Y)$ the right adjoint to the inverse image $f^{-1}:{\cal O}(Y)\to {\cal O}(X)$; notice that, concretely, $f_{\ast}$ sends any open set $U$ of $X$ to the smallest open set of $Y$ containing the image $f(U)$ of $U$ under $f$.

\begin{corollary}
Let $X$ be a (completely regular) topological space and $\eta_{X}:X\to \beta(X)$ be its Stone-\v{C}ech compactification; let $B$ be a normal Wallman base for $X$ such that the sets of the form ${\eta_{X}}_{\ast}(U)$ where $U\in B$ form a sublattice of $\beta(X)$ which is a base for it (e.g. $B$ is a normal Wallman base of $X$ which contains $Coz(X)$). Then, under the assumption that the maximal ideal theorem holds, $\beta(X)$ is homeomorphic to $Max(B)$. 
\end{corollary}

\begin{proofs}
Under the hypotheses of the theorem, ${\cal O}(X)$ is a completely regular space such that $B$ is a normal ${\cal O}(X)$-conjunctive base for it (cf. the remarks following the definition of $A$-conjunctive sublattice); hence the hypotheses of Theorem \ref{genwallman} are satisfied and we obtain an isomorphism $\beta({\cal O}(X))\simeq Id_{J_{m}^{B}(B)}$. Passing to the spaces of points of these locales we obtained our desired homeomorphism $X\cong Max(B)$. 
\end{proofs}  

If we identify $\beta(X)$, as in Corollary IV 2.3 \cite{stone}, with the space of maximal completely regular filters on ${\cal O}(X)$ then the map $\eta_{X}:X\to \beta(X)$ can be described as the function sending any point $x\in X$ to the filter $\{U\in {\cal O}(X) \textrm{ | } x\in U\}$. By the description of $\beta(X)$ given in section IV 2.3 \cite{stone}, one immediately deduces that $\eta'$ sends a point $x\in X$ to the collection of opens of the form $\bigcup I$ where $I$ is a completely regular ideal such that there exists $v\in I$ such that $x\in v$. Now, it is clear that every such open contains $x$; conversely, given $v\in {\cal O}(X)$ such that $x\in V$, consider the completely regular ideal $I_{v}:=\{u\in {\cal O}(X) \textrm{ | } U \eqslantless v\}$, where $\eqslantless$ denotes the completely below relation on ${\cal O}(X)$; $X$ being completely regular, $v=\bigcup I_{v}$ and hence, since $x\in v$, there exists $z\in I_{v}$ such that $x\in z$. 
For any open set $U$ of $X$, ${\eta_{X}}_{\ast}(U)$ can be identified with the open set $\{M\in {\cal M}_{{\cal O}(X)} \textrm{ | } U\in M\}$ and ${\eta_{X}}_{\ast}$ defines an injective distributive lattice homomorphism ${\cal O}(X)\to {\cal O}(\beta(X))$.
The homeomorphism $\beta(X)\to Max(B)$ provided by the Corollary can be identified with the map sending a maximal completely regular filter on ${\cal O}(X)$ to the complement in $B$ of its intersection with $B$.

Let us now interpret Theorem \ref{genwallman} in the context of specific constructions of the Stone-\v{C}ech compactification of a locale. Recall from \cite{stone} (Theorem IV 2.2) that the Stone-\v{C}ech compactification of a locale $A$ can be identified with the locale map $\eta_{A}:A\to C(A)$, where $C(A)$ is the locale of completely regular ideals of $A$ and the direct image ${\eta_{A}}_{\ast}$ of $\eta_{A}$ sends any element $a\in A$ to the smallest completely regular ideal of $A$ containing the principal ideal $(a)$. Let us suppose that $A$ and $B$ satisfy the hypotheses of Theorem \ref{genwallman}. Then we have a locale isomorphism $h:C(A)\to Id_{J_{m}^{B}}(B)$ such that $h\circ \eta_{A}=g$, where $g:A\to Id_{J_{m}^{B}}(B)$ is the locale map whose direct image sends any element $a\in A$ to the $J_{m}^{B}$-closure of the ideal $\{b\in B \textrm{ | } b\leq a\}$. Since $h$ is an isomorphism then $h_{\ast}$ preserves arbitrary joins, whence $h_{\ast}$ can be identified with the map sending any completely regular ideal $I$ on $A$ to the $J_{m}^{B}$-closure of the ideal $\{b\in B \textrm{ | } b\leq a \textrm{ for some } a\in I\}$. We may thus conclude that this map is bijective. Notice that in the case of a normal subfit locale $A$ such map is the function sending to every completely regular ideal $I$ of $A$ to itself, regarded as an element of $Id_{J_{m}^{B}}(B)$; therefore for any normal subfit locale $A$, the completely regular ideals of $A$ are precisely the $J_{m}^{A}$-closed ones. Specifically, we have the following result.

\begin{corollary}
Under the hypotheses of Theorem \ref{genwallman}, the map $C(A)\to Id_{J_{m}^{B}}(B)$ sending any completely regular ideal $I$ $A$ to the $J_{m}^{B}$-closure of the ideal $\{b\in B \textrm{ | } b\leq a \textrm{ for some } a\in I\}$ is a bijection.  

In particular, for any normal subfit locale $A$ and any ideal $I$ of $A$, $I$ is completely regular if and only if for any $b\in A$ such that for any $c\in A$, $c\vee b=1$ implies $c\vee d=1$ for some $d\in I$, $b\in I$.
\end{corollary}\qed  

An alternative description of the Stone-\v{C}ech compactification of a locale is given in \cite{BMSC3}; specifically, for any locale $A$ its Stone-\v{C}ech compactification $\beta(A)$ is characterized as the geometric syntactic category of a geometric propositional theory which axiomatizes the maximal completely regular filters on $A$. In particular, in presence of the axiom of choice such locale can be described as the topological space whose underlying set is the collection ${\cal M}_{A}$ of all maximal completely regular filters on $A$ and whose topology is generated by basic open sets of the form ${\cal M}_{A}^{a}:=\{M\in {\cal M}_{A} \textrm{ | } a\in A\}$ for $a\in A$. The direct image of the canonical map $\eta_{A}:A\to \beta(A)$ can be identified with the function sending any $a\in A$ to the open set ${\cal M}_{A}^{a}$. It is easy to see that, under the hypotheses of Theorem \ref{genwallman}, the homeomorphism $A\cong Max(B)$ can be identified with the map sending any maximal completely regular filter on $A$ to the complement in $B$ of its intersection with $B$. In particular, we have the following result.

\begin{corollary}
Under the hypotheses of Theorem \ref{genwallman}, the map ${\cal M}_{A}\to Id(B)$ sending any maximal completely regular filter $F$ on $A$ to $B\setminus (F\cap B)$ is a bijection onto $Max(B)$.  
\end{corollary} 

In passing, we note that the space $\beta(A)$ admits the following topos-theoretic description. For any locale $A$, one can define a Grothendieck topology $J_{A}^{cr}$ on $A$, called the \emph{completely regular topology} on $A$, as follows: for any $a\in A$ and any sieve $S$ on $a$, $S\in J_{A}^{cr}(a)$ if and only if for every $b\in A$ such that $b \eqslantless A$ in $A$, the arrow $b\leq a$ belongs to $S$ (here $\eqslantless$ denotes the completely below or `really inside' relation on $A$ of section IV 1.4 of \cite{stone}). Let us verify that this is indeed a Grothendieck topology: the maximality axiom is clearly satisfied, the stability axiom holds since for any $a, b\in A$ with $b\leq a$ for any $c\in A$ such that $c \eqslantless b$, $c \eqslantless a$ and hence $c=c\wedge b \leq b$ belongs to the pullback along $b\leq a$ of any $J_{A}^{cr}$-covering sieve on $a$, and the transitivity axiom holds as a consequence of the subdivisibility property of the relation $\eqslantless$ (cf. Lemma IV 1.4 \cite{stone} for the properties of the relation $\eqslantless$). 

Notice that for any locale $A$, the points of the topos $\Sh(A, J_{A}^{cr})$ are exactly the completely regular filters of $A$ (in the sense of \cite{stone}). Observing that the specialization order on the space of points $X_{A}^{cr}$ of $\Sh(A, J_{A}^{cr})$ can be identified with the subset-inclusion relation between completely regular filters of $A$, we conclude that the maximal completely regular filters of $A$ can be characterized as the points of the $\Sh(A, J_{A}^{cr})$ which are maximal with respect to the specialization ordering on $X_{A}^{cr}$, that is as the points $x$ of $X_{A}^{cr}$ such that the intersection of all the open sets of $X_{A}^{cr}$ containing $x$ is equal to $\{x\}$. Let us denote by $X_{A}^{mcr}$ the subspace of $X_{A}^{cr}$ consisting of such points; then $X_{A}^{cr}$ is the topological space whose underlying set is the set ${\cal M}_{A}$ of maximal completely regular filters of $A$ and whose topology is generated by the basic open sets of the form ${\cal M}_{A}^{a}:=\{M\in {\cal M}_{A} \textrm{ | } a\in A\}$ for $a\in A$, and hence can be identified with the Stone-\v{C}ech compactification of $A$. Notice that we have a commutative diagram
    
\[  
\xymatrix {
\Sh(A, J_{A}^{can}) \ar[dr]_{\Sh(\eta_{A})} \ar[r]^{f} & \Sh(A, J_{A}^{cr}) \\
 & \Sh(X_{A}^{cr}) \ar[u]^{i},}
\] 
where $i:\Sh(X_{A}^{mcr})\hookrightarrow \Sh(X_{A}^{cr})$ is the geometric inclusion induced by the subspace inclusion $X_{A}^{mcr}\subseteq X_{A}^{cr}$ and $f:\Sh(A, J_{A}^{can}) \to \Sh(A, J_{A}^{cr})$ is the geometric morphism corresponding to the meet-semilattice homomorphism $A\to A$ sending any element $a\in A$ to the join $\mathbin{\mathop{\textrm{\huge $\vee$}}\limits_{b \eqslantless a}}b$. Notice that if $A$ is completely regular then the topology $J_{A}^{cr}$ is subcanonical and the geometric morphism $f$ can be identified with the morphism induced by the morphism of sites $(A, J_{A}^{cr}) \to (A, J_{A}^{can})$ given by the identity map on $A$; in particular, $f$ is a geometric inclusion.  

In view of the general techniques introduced in \cite{OC10}, it is worth to remark that all these different constructions of the Stone-\v{C}ech compactification of a locale admit a natural topos-theoretic interpretation as Morita-equivalences between different geometric theories (equivalently, as different site representations for the same topos), which can be fruitfully exploited to systematically transfer properties and results between one presentation and another and hence to effectively tackle questions about the Stone-\v{C}ech compactification of a locale from several different points of view (cf. \cite{OC11} for a comprehensive analysis of the importance of Morita-equivalences and an overview of the general methodologies for transferring results across them). Specifically, for any locale $A$, we have a representation
\[
\Sh(\beta(A))\simeq \Sh(C(A), J_{C(A)}^{can}),
\]
corresponding to the construction of $\beta(A)$ as the locale $C(A)$. 

Also, we have a Morita-equivalence
\[
\Sh(\beta(A))\simeq \Set[Max{\mathbb R}(A)],
\]
corresponding to the construction of $\beta(A)$ as the locale of maximal ideals of the ring of bounded locale maps from $A$ to the locale $\mathbb R$ of real numbers, where $\Set[Max{\mathbb R}(A)]$ is the classifying topos of the propositional geometric theory $Max{\mathbb R}(A)$ axiomatizing such ideals introduced in \cite{BMSC2}.

Alternatively, the Stone-\v{C}ech compactification of a locale $A$ can be built as the canonical locale map from $A$ to the locale given by the geometric syntatic category of the propositional theory ${\mathbb M}_{A}$ of almost prime (equivalently, maximal) completely regular filters on $A$ defined in \cite{BMSC2}. This representation can be expressed as an equivalence
\[
\Sh(\beta(A))\simeq \Set[{\mathbb M}_{A}]
\]    
between $\Sh(\beta(A))$ and the classifying topos $\Set[{\mathbb M}_{A}]$ of this theory. 

Our construction of the Stone-\v{C}ech compactification of $A$ as a subtopos of $\Sh(A, J_{A}^{cr})$ established above yields an additional Morita-equivalence
\[
\Sh(\beta(A))\simeq \Sh(X_{A}^{cr}).
\]

Also, if $B$ is a base of $A$ satisfying the hypotheses of Theorem \ref{genwallman}, we have a further representation
\[
\Sh(\beta(A))\simeq \Sh(B, J_{B}^{m}).
\]

\section{Dualities between topological spaces and distributive lattices}\label{dualtopoldlat}

\subsection{The general framework}\label{generalframework}

Let us define $\textbf{TopDLat}$ as the category whose objects are the pairs $(X, D)$ where $X$ is a topological space and $D$ is a sublattice of ${\cal O}(X)$ and whose arrows $(X, D)\to (Y, D')$ are the continuous maps $f:X\to Y$ such that the inverse image $f^{-1}:{\cal O}(Y)\to {\cal O}(X)$ restricts to (a distributive lattice homomorphism) $D'\to D$:
\[  
\xymatrix {
{\cal O}(Y) \ar[r]^{f^{-1}} & {\cal O}(X)\\
D' \ar[u] \ar[r]^{f^{-1}|D'} & D \ar[u]}
\]
where the arrows $D'\to {\cal O}(Y)$ and $D\to {\cal O}(X)$ are the canonical inclusions.

Let us denote by $\textbf{TopDLat}_{W}$ the full subcategory of $\textbf{TopDLat}$ on the objects $(X, B)$ such that $X$ is a $T_{0}$-space, $B$ is a Wallman base for $X$ and the map $\eta^{X}_{B}:X\to Spec(B)$ is surjective (equivalently, a homeomorphism) onto $Max(B)$. 

Let us define $\textbf{DLat}_{W}$ as the subcategory of the category $\textbf{DLat}$ of distributive lattices and distributive lattice homomorphisms between them whose objects are the conjunctive distributive lattices $D$ and whose arrows are the maximal homomorphisms between them (cf. section \ref{maxspectrum} above for the definition of maximal homomorphism). 

We can define two functors
\[
H:\textbf{TopDLat}_{W}^{\textrm{op}} \to \textbf{DLat}_{W}
\]  
and
\[
K:\textbf{DLat}_{W} \to \textbf{TopDLat}_{W}^{\textrm{op}} 
\]  
as follows. 

For any $(X, B)\in \textbf{TopDLat}_{W}$ we set $H((X, B))=B$ and for any arrow $f:(X, B)\to (Y, B')$ in $\textbf{TopDLat}_{W}$ we set $H(f)$ equal to the restriction $f^{-1}:B'\to B$ of the inverse image map $f^{-1}:{\cal O}(Y)\to {\cal O}(X)$. 

Let us verify that $H$ is well-defined. We have to show that for any $(X, B)\in \textbf{TopDLat}_{W}$, $B$ is conjunctive, that is the maximal topology $J_{m}^{B}$ on $B$ is subcanonical. We observe that 
\[
\Sh(\eta^{X}_{B}):\Sh(X)\to \Sh(Spec(B))
\]
is isomorphic (as a subtopos of $\Sh(Spec(B)$) to the geometric morphism
\[
\Sh(B, J^{can}_{{\cal O}(X)}|B) \to \Sh(B, J^{coh}_{B})\simeq \Sh(Spec(B)) 
\]
induced by the morphism of sites $(B, J^{coh}_{B})\to (B, J^{can}_{{\cal O}(X)}|B)$.

From this it follows that if $\eta^{X}_{B}:X \to Max(B)$ is an homeomorphism then the geometric morphism $\Sh(\eta^{X}_{B}):\Sh(X)\to \Sh(Spec(B))$ is isomorphic to the geometric morphism $\Sh(Max(B))\to \Sh(Spec(B))$ induced by the subspace inclusion $Max(B)\hookrightarrow Spec(B)$ and hence to the geometric morphism
\[
\Sh(B, J^{B}_{m}) \to \Sh(B, J^{coh}_{B})\simeq \Sh(Spec(B)) 
\]
induced by the morphism of sites $(B, J^{coh}_{B})\to (B, J^{B}_{m})$.

We can thus conclude, from the isomorphism of subtoposes
\[
\Sh(B, J^{can}_{{\cal O}(X)}|B) \to \Sh(B, J^{coh}_{B})\simeq \Sh(Spec(B)) 
\]    
and 
\[
\Sh(B, J^{B}_{m}) \to \Sh(B, J^{coh}_{B})\simeq \Sh(Spec(B)) 
\]
the equality of Grothendieck topologies $J_{m}^{B}=J^{can}_{{\cal O}(X)}|B$, from which it follows in particular that $J_{m}^{B}$ is subcanonical, i.e. that $B$ is conjunctive, as required. 

From the definition of the arrows in the category $\textbf{TopDLat}_{W}$ and the equalities $J_{m}^{B}=J^{can}_{{\cal O}(X)}|B$ just observed it immediately follows that for any arrow $f:(X, B)\to (Y, B')$ in $\textbf{TopDLat}_{W}$ the restriction $f^{-1}|_{B'}:B'\to B$ is a maximal homomorphism, that is an arrow $B'\to B$ in $\textbf{DLat}_{W}$. 

We can thus conclude that the functor $H$ indeed takes values in the category $\textbf{DLat}_{W}$.

Let us now show that the functor $K:\textbf{DLat}_{W} \to \textbf{TopDLat}_{W}^{\textrm{op}}$ is well-defined.

We have to show that for any $D$ in $\textbf{DLat}_{W}$, $(Max(D), D)$ is an object of $\textbf{TopDLat}_{W}^{\textrm{op}}$. Since $Max(D)$ is a $T_{1}$-space (cf. Lemma II3.5 \cite{stone}) then it is in particular a $T_{0}$-space; therefore $(Max(D), D)$ is an object of $\textbf{TopDLat}_{W}^{\textrm{op}}$, as required. Given an arrow $f:D\to D'$ in $\textbf{DLat}_{W}$, that is a maximal homomorphism $D\to D'$, we have to show that the induced map $Max(f):Max(D')\to Max(D)$ yields an arrow $(Max(D'), D')\to (Max(D), D)$ in the category $\textbf{TopDLat}_{W}$. To this end, we observe that the injections $s_{D}:D\mono {\cal O}(Max(D))$ and $s_{D'}:D'\mono {\cal O}(Max(D'))$) respectively given by the composite of the canonical injections $D\mono Id_{J_{m}^{D}}(D)$ and $D'\mono Id_{J_{m}^{D'}}(D')$ with the isomorphisms $Id_{J_{m}^{D}}(D) \cong {\cal O}(Max(D))$ and $Id_{J_{m}^{D'}}(D') \cong {\cal O}(Max(D'))$ make the following diagram commute:    
  
\[  
\xymatrix {
{\cal O}(Max(D')) \ar[r]^{Max(f)^{-1}} & {\cal O}(Max(D))\\
D' \ar[u]^{s_{D'}} \ar[r]^{f} & D \ar[u]^{s_{D}}.}
\]  
  
\begin{theorem}\label{duality}
With the notation above, the functors
\[
H:\textbf{TopDLat}_{W}^{\textrm{op}} \to \textbf{DLat}_{W}
\]  
and
\[
K:\textbf{DLat}_{W} \to \textbf{TopDLat}_{W}^{\textrm{op}} 
\]
are inverse to each other (up to natural isomorphism) and hence yield a duality between the category  $\textbf{TopDLat}_{W}$ and the category $\textbf{DLat}_{W}$.    
\end{theorem} 

\begin{proofs}
For any arrow $(X, D)\to (Y, D')$ in $\textbf{TopDLat}$ we have a commutative diagram
\[  
\xymatrix {
X \ar[d]^{f} \ar[r]^{\eta^{X}_{D}} & Spec(D) \ar[d]^{Spec(j)}\\
Y \ar[r]_{\eta^{Y}_{D'}} & Spec(D')}
\]  
where $j:D'\to D$ is the restriction $D'\to D$ of $f^{-1}:{\cal O}(Y)\to {\cal O}(X)$. Indeed, for any prime ideal $P$ of $D$, $Spec(j)(P)=\{v\in {\cal O}(Y) \textrm{ | } f^{-1}(v)\in P\}$ and hence for any $x\in X$, $Spec(j)(\eta^{X}_{D}(x))=\eta^{Y}_{D'}(f(x))$. 

By the commutativity of this square, the map $\eta^{X}_{D}:X\to Max(D)$ (for $(X, D)\in \textbf{TopDLat}_{W}^{\textrm{op}}$,) is an isomorphism in $\textbf{TopDLat}_{W}$ which is natural in $(X, D)\in \textbf{TopDLat}_{W}^{\textrm{op}}$; hence it gives rise to a natural isomorphism $1_{\textbf{TopDLat}_{W}^{\textrm{op}}}\to K\circ H$. Conversely, for any $D$ in $\textbf{TopDLat}_{W}$, the isomorphism identifying $D$ with its image in ${\cal O}(Max(D))$ under the embedding $D\mono {\cal O}(Max(D))$ is natural in $D\in \textbf{TopDLat}_{W}$ and hence yields a natural isomorphism $1_{\textbf{TopDLat}_{W}}\to H\circ K$.      
\end{proofs}

We might naturally wonder whether we can naturally characterize the distributive lattices $D$ which correspond to compact Hausdorff spaces under this duality, that is the distributive lattices $D$ such that the space $Max(D)$ is Hausdorff. In $\cite{almax}$ Johnstone gave the following characterization: $Max(D)$ is Hausdorff if and only if $D$ is \emph{semi-normal}, i.e. for any $a,b\in D$ such that $a\vee b=1$ there exist $c,d\in D$ such that $c\vee a=b\vee d=1$ and $c\wedge d$ belongs to the $J_{m}^{D}$-closure $\overline{(0)}^{J_{m}^{D}}$ of the principal ideal $(0)$. Notice that for any $a\in D$ the $J_{m}^{D}$-closure $\overline{(a)}^{J_{m}^{D}}$ of a principal ideal $(a)$ is the ideal $\{c\in D \textrm{ | for all $b$ such that } b\vee c=1, b\vee a=1 \}$; therefore the condition $c\wedge d\in \overline{(0)}^{J_{m}^{D}}$ can be reformulated by saying that for all $e\in D$ with $e\vee (c\wedge d)=1$, $e=1$. 

In passing, we record the following fact.

\begin{proposition}
For any conjunctive distributive lattice $D$ (for example, a Wallman base for some topological space), $D$ is normal if and only if it is semi-normal.
\end{proposition}\qed
 
A common situation in which Theorem \ref{duality} can be applied is when we have a category $\cal K$ of topological spaces $X$, each of which equipped with a Wallman base $B_{X}$ with the property that $\eta^{X}_{B_{X}}:X \to Max(B_{X})$ is an homeomorphism and for any arrow $f:X\to Y$ in $\cal K$ the inverse image $f^{-1}:{\cal O}(Y)\to {\cal O}(X)$ restricts to a map $B_{Y}\to B_{X}$. Indeed, in such a situation one has a functor $Z:{\cal K}\to \textbf{TopDLat}_{W}$, sending any space $X$ in $\cal K$ to the pair $(X, B_{X})$ and any arrow $f:X\to Y$ in $\cal K$ to the map $f:(X, B_{X})\to (Y, B_{Y})$, which identifies $\cal K$ with a full subcategory of $\textbf{TopDLat}_{W}$. Hence the restriction $H|{\cal K}^{\textrm{op}}:{\cal K}^{\textrm{op}}\to \textbf{DLat}_{W}$ of $H$ to ${\cal K}^{\textrm{op}}$ will yield an equivalence between ${\cal K}^{\textrm{op}}$ and the full subcategory $ExtIm(H)$ of $\textbf{DLat}_{W}$ on the objects which are isomorphic to one of the form $H(X, B_{X})$ (for $X\in {\cal K}$).

We shall see a couple of instances of this phenomenon in the next two sections.

\subsection{A duality for compact Hausdorff spaces}\label{dualityalexandrovalg}

For any compact Hausdorff space $X$ the sublattice $Coz(X)$ of ${\cal O}(X)$ consisting of the sets of the form $Coz(f)=f^{-1}(\mathbb R\setminus \{0\})$ for a continuous function $f:X\to {\mathbb R}$ $Coz(X)$ is a normal Wallman base for $X$ (cf. Propositions IV2.6 and IV3.3 \cite{stone}) such that $X\cong Max(Coz(X))$; moreover, any continuous map of topological spaces $f:X\to Y$ induces a distributive lattice homomorphism $f^{-1}:Coz(Y)\to Coz(X)$. Therefore the method of the last section yields a duality between the category $\textbf{CHaus}$ of compact Hausdorff spaces and continuous maps between them and the full subcategory $\textbf{CHDLat}$ of $\textbf{DLat}_{W}$ on the distributive lattices $D$ which are, up to isomorphism, of the form $Coz(X)$ for some compact Hausdorff space $X$. We can characterize these distributive lattices more intrinsically as follows. For a compact Hausdorff space $X$, the inclusion $Coz(X)\hookrightarrow {\cal O}(X)$ corresponds, under the isomorphism ${\cal O}(X)\cong {\cal O}(Max(Coz(X)))$, to the morphism $Coz(X)\to {\cal O}(Max(Coz(X)))$ sending any element $d\in Coz(X)$ to the open set $\{M\in Max(Coz(X)) \textrm{ | } d\notin M\}$. We can thus conclude that a distributive lattice $D$ in $\textbf{DLat}_{W}$ is isomorphic to one of the form $Coz(X)$ if and only if the image of the canonical map $D \mono {\cal O}(Max(D))$ coincides with $Coz(Max(D))$, that is if and only if for any $d\in D$ there exists a continuous map $f:Max(D)\to {\mathbb R}$ such that for any $M\in Max(D)$, $f(M)=0$ if and only if $d\in M$, and for any continuous map $g:Max(D)\to {\mathbb R}$ there exists a (necessarily unique) element $d\in D$ such that for any $M\in Max(D)$, $g(M)=0$ if and only if $d\in M$. 

Summarizing, we have the following result.

\begin{theorem}\label{dualdlat}
The functors 
\[
Coz:\textbf{CHaus}^{\textrm{op}}\to \textbf{CHDLat}
\]
and 
\[
Max:\textbf{CHDLat} \to \textbf{CHaus}^{\textrm{op}} 
\]
defined above yield a duality between the category $\textbf{CHaus}$ of compact Hausdorff spaces and the subcategory $\textbf{CHDLat}$ of the category of distributive lattices.
\end{theorem}\qed

We can characterize the category $\textbf{CHDLat}$ more explicitly thanks to the notion of Alexandrov algebra. Recall from \cite{stone} that an \emph{Alexandrov algebra} is a normal distributive lattice $D$ in which countable joins exist and distribute over finite meets and the following `approximation property' holds: for any $a \in D$, there exist sequences $\{b_{n} \textrm{ | } n\in {\mathbb N}\}$ and $\{c_{n} \textrm{ | } n\in {\mathbb N}\}$ of elements of $D$ such that $\mathbin{\mathop{\textrm{\huge $\vee$}}\limits_{n\in {\mathbb N}}}c_{n}=a$, $b_{n}\wedge c_{n}=0$ and $b_{n}\vee a=1$ for all $n\in {\mathbb N}$.

On any Alexandrov algebra $D$ one can define a Grothendieck topology $C_{D}$, called the \emph{countable topology} on $D$, by saying that the $C_{D}$-covering sieves on a given element $d$ are precisely the sieves which contain countable families of arrows whose join is equal to $d$. We define the category $\textbf{AlexAlg}$ of Alexandrov algebras as the category whose objects are the Alexandrov algebras and whose morphisms are the distributive lattice homomorphisms between them which preserve countable joins.

It is proved in \cite{stone} (cf. also \cite{reynolds1} and \cite{reynolds2} for the original sources) that for any Alexandrov algebra $D$, the frame $Id_{C_{D}}(D)$ of $C_{D}$-ideals on $D$ is a completely regular locale. 

We shall now prove that the lattices in $\textbf{CHDLat}$ can be identified with the Alexandrov algebras $D$ such that the frame $Id_{C_{D}}(D)$ is compact, while the morphisms in $\textbf{CHDLat}$ coincide precisely with the morphisms of Alexandrov algebras (that is, with the arrows in $\textbf{AlexAlg}$). This will provide us with a more intrinsic lattice-theoretic duality for compact Hausdorff spaces.  

To this end, we prove the following result.

\begin{lemma}
Let $({\cal C}, J)$ be a site and $c$ be an object of $\cal C$. Then the $J$-closure $\overline{(c)}^{J}$ of the principal ideal $(c):=\{d\in {\cal C} \textrm{ | there is } f:d\to c \textrm{ in $\cal C$} \}$ in $\cal C$ generated by $c$ satisfies the property that every covering of $\overline{(c)}^{J}$ by subterminals in $\Sh({\cal C}, J)$ admits a finite subcovering if and only if any $J$-covering sieve on $c$ contains a $J$-covering sieve generated by a finite family of arrows.
\end{lemma}

\begin{proofs}
Since every subterminal of $\Sh({\cal C}, J)$ can be written as a join of subterminals of the form $\overline{(a)}^{J}$ for some $a\in {\cal C}$, every covering of $\overline{(c)}^{J}$ by subterminals in $\Sh({\cal C}, J)$ admits a finite subcovering if and only if every covering of $\overline{(c)}^{J}$ by subterminals of the form $\overline{(a)}^{J}$ admits a finite subcovering. Now, given a family $\{c_{i} \textrm{ | } i\in I\}$ of objects of $\cal C$, $\overline{(c)}^{J}=\mathbin{\mathop{\textrm{\huge $\vee$}}\limits_{i\in I}}\overline{(c_{i})}^{J}$ in $\Sh({\cal C}, J)$ if and only if there exists a $J$-covering sieve $S$ on $c$ such that for any $f\in S$, $dom(f)\in (c_{i})$ for some $i\in I$, equivalently the sieve $\{f:dom(f)\to c \textrm{ | } dom(f)\in (c_{i}) \textrm{ for some $i\in I$}\}$ is $J$-covering on $c$. From this our thesis immediately follows.  
\end{proofs}

Applying the lemma to the site $(D, C_{D})$ where $D$ is an Alexandrov algebra and $C_{D}$ is the countable topology on it (with $c$ equal to the top element $1$ of $D$) yields the following criterion: the frame $Id_{C_{D}}(D)$ is compact if and only if every $C_{D}$-covering sieve on $1$ contains a $C_{D}$-covering sieve generated by a finite family of arrows; in other words, $Id_{C_{D}}(D)$ is compact if and only if for any denumerable family $\{c_{n} \textrm{ | } n\in {\mathbb N}\}$ of elements in $D$ such that $\mathbin{\mathop{\textrm{\huge $\vee$}}\limits_{n\in {\mathbb N}}} c_{n}=1$ there exists a finite subset $H\subseteq {\mathbb N}$ such that $\mathbin{\mathop{\textrm{\huge $\vee$}}\limits_{n\in H}} c_{n}=1$. We shall say that an Alexandrov algebra is \emph{countably compact} if it satisfies this condition. 

Notice that, for any Alexandrov algebra $D$, since the frame $Id_{C_{D}}(D)$ is completely regular, the space $X_{(D, C_{D})}$ of its points is Hausdorff (cf. III1.1 \cite{stone}). Concretely, this means that for any two distinct $C_{D}$-prime filters $F$ and $G$ on $D$ there exist elements $c,d\in D$ such that $c\in F$, $d\in G$ and $c\wedge d=0$, or equivalently that the subsets of the form $\{F \in X_{(D, C_{D})} \textrm{ | } c\in F\}$ (for $c\in D$) form a (normal) Wallman base for $X_{(D, C_{D})}$ (i.e., for any $c\in D$ and any $C_{D}$-prime filter $F$ on $D$ such that $c\in F$ there exists $d\in D$ such that $c\vee d=1$ and $d\notin F$).

Next, we proceed to show that the category $\textbf{CHDLat}$ coincides with the full subcategory $\textbf{AlexAlg}_{c}$ of $\textbf{AlexAlg}$ on the countably compact Alexandrov algebras. Let us start by proving that $\textbf{CHDLat}$ is a full subcategory of $\textbf{AlexAlg}$. Let us thus suppose $X$ to be a compact Hausdorff space; we want to show that $Coz(X)$ is a countably compact Alexandrov algebra. The fact that $Coz(X)$ is an Alexandrov algebra was proved in IV2.9 \cite{stone} in the more general case of a completely regular space $X$; it thus remains to show that $Coz(X)$ is countably compact. But, since $Coz(X)$ is closed in ${\cal O}(X)$ under countable unions (cf. Lemma IV2.5 \cite{stone}), the fact that $X$ is compact implies that $Coz(X)$ is countably compact, since $X$ is the top element of $Coz(X)$.      

Given a morphism $f:D\to D'$ in $\textbf{CHDLat}$, let us show that it is an arrow $D\to D'$ in $\textbf{AlexAlg}$, that is a morphism of sites $(D, C_{D})\to (D', C_{D'})$. The morphism $f$ induces a continuous map $Max(f):Max(D')\to Max(D)$ such that the restriction $f^{-1}|:Coz(Max(D')) \to Coz(Max(D))$ of its inverse image is isomorphic to $f$; clearly $f^{-1}$ preserves countable (in fact, arbitrary) unions in ${\cal O}(Max(D'))$ and since $Coz(Max(D'))$ (resp. $Coz(Max(D))$) is closed in ${\cal O}(Max(D'))$ (resp. in ${\cal O}(Max(D))$) under countable unions (cf. Lemma IV2.5 \cite{stone}) $f^{-1}|\cong f$ preserves countable joins and therefore is an arrow $D\to D'$ in $\textbf{AlexAlg}$.

So far we have proved that the category $\textbf{CHDLat}$ is a subcategory of the category $\textbf{AlexAlg}_{c}$. It therefore remains to prove the converse. Given a countably compact Alexandrov algebra $D$, as we saw above $D$ is isomorphic to the co-zero set $Coz(X_{(D, C_{D})})$ of the compact Hausdorff space $X_{(D, C_{D})}$ and hence by Theorem \ref{dualdlat} it is an object of the category $\textbf{CHDLat}$; clearly, this holds functorially in $D\in \textbf{AlexAlg}_{c}$ (by Proposition IV2.10 \cite{stone}), from which it follows that the category $\textbf{AlexAlg}_{c}$ is a subcategory of $\textbf{CHDLat}$. 

Therefore we can conclude that the categories $\textbf{CHDLat}$ and $\textbf{AlexAlg}_{c}$ are equal.   

The assignment $D\to X_{(D, C_{D})}$ can be made into a functor 
\[
C:\textbf{AlexAlg}_{c} \to \textbf{CHaus}^{\textrm{op}};
\]
indeed, any arrow $f:D\to D'$ in $\textbf{AlexAlg}_{c}$ induces a locale morphism $Id_{C_{D'}}(D')\to Id_{C_{D'}}(D')$ and hence a continuous map $C(f):X_{(D', C_{D'})} \to X_{(D, C_{D})}$ between the spaces of points of the two locales. 

The functor $Max:\textbf{CHDLat} \to \textbf{CHaus}^{\textrm{op}}$ is naturally isomorphic to the functor $C:\textbf{AlexAlg}_{c} \to \textbf{CHaus}^{\textrm{op}}$. Indeed, for any countably compact Alexandrov algebra $D$, the space $X_{(D, C_{D})}$ is compact Hausdorff with normal Wallman base $Coz(X_{(D, C_{D})}) \cong D$ and hence it is homeomorphic to $Max(D)$ (cf. Corollary \ref{compWallman}), naturally in $D$.

Summarizing, we have the following duality theorem.

\begin{theorem}\label{dualalex}
The functors 
\[
Coz:\textbf{CHaus}^{\textrm{op}}\to \textbf{AlexAlg}_{c}
\]
and 
\[
C:\textbf{AlexAlg}_{c} \to \textbf{CHaus}^{\textrm{op}} 
\]
defined above yield a duality between the category $\textbf{CHaus}$ of compact Hausdorff spaces and the full subcategory $\textbf{AlexAlg}_{c}$ of the category of Alexandrov algebras on the countably compact ones. In fact, the functor 
\[
C:\textbf{AlexAlg}_{c} \to \textbf{CHaus}^{\textrm{op}} 
\]
is naturally isomorphic, under the identification $\textbf{AlexAlg}_{c}=\textbf{CHDLat}$, to the functor
\[
Max:\textbf{CHDLat} \to \textbf{CHaus}^{\textrm{op}} 
\]
of Theorem \ref{dualdlat} above.
\end{theorem}

Note that for any topological space $X$, $Coz(X)$ is an Alexandrov algebra (cf. \cite{stone}) and if $X$ is compact then $Coz(X)$ is countably compact (cf. the argument given above). 

As an immediate corollary of our duality theorem we obtain the following result about Alexandrov algebras.

\begin{proposition}
Let $D$ be an Alexandrov algebra. Then $D$ is conjunctive, with the maximal topology $J_{m}^{D}$ equal to the countable topology $C_{D}$ (i.e., for any sieve $S:=\{c_{i}\leq c \textrm{ | } \in I\}$ on $c\in D$ with the property that for every $d\in D$, $c\vee d=1$ implies that there exists a finite subset $J\subseteq I$ such that $\mathbin{\mathop{\textrm{\huge $\vee$}}\limits_{i\in J}}c_{i} \vee d=1$, we have $\mathbin{\mathop{\textrm{\huge $\vee$}}\limits_{i\in K}}c_{i}=c$ for some countable subset $K\subseteq I$), if and only if it is countably compact.  
\end{proposition}\qed    

In light of Gelfand duality between commutative $C^{\ast}$-algebras and compact Hausdorff spaces, the duality of Theorem \ref{dualalex} also provides an equivalence between the category of $C^{\ast}$-algebras and the category of countably compact Alexandrov algebras, obtained by composition of the two dualities. We shall now give an explicit description of this categorical equivalence which allows to (functorially) construct the Alexandrov algebra associated to a given commutative $C^{\ast}$-algebra directly in terms of it, and conversely to explicitly construct the $C^{\ast}$-algebra corresponding to a given countably compact Alexandrov algebra in terms of it. 

We shall describe this equivalence in the case of \emph{real} $C^{\ast}$-algebras (that is, of rings of real-valued continuous functions on a compact Hausdorff space, cf. \cite{Stone} and IV 4.4 \cite{stone} for an axiomatic description of this notion), but our arguments can be straightforwardly extended to the context of complex $C^{\ast}$-algebras to yield an equivalence between the usual category of (complex) $C^{\ast}$-algebras and the category $\textbf{AlexAlg}_{c}$. 

Given a $C^{\ast}$-algebra $A$, by Gelfand duality we have a canonical surjective lattice homomorphism $\eta_{L(A)}:L(A)\to Coz(Max(L(A)))\cong Coz(Max(A))$, which induces an equivalence 
\[
\Sh(Max(L(A)))\simeq \Sh(Max(A))\simeq \Sh(Coz(Max(A))).
\]
Therefore $Coz(Max(A))$ can be characterized up to isomorphism as the image of the map $\eta_{L(A)}$. More precisely, by Gelfand duality, the map $A\to {\cal O}(Max(A))$ sending any $a\in A$ to the open set $\{M\in Max(A) \textrm{ | } a\notin M\}$ is surjective and factors through the canonical map $A\to L(A)$ yielding the morphism $\eta_{L(A)}$; hence $Coz(Max(A))$ can be identified with the sublattice of ${\cal O}(Max(A))$ consisting of the subsets of the form $\{M\in Max(A) \textrm{ | } a\notin M\}$ (for $a\in A$). This sublattice is therefore an Alexandrov algebra $D_{A}$, which we call the Alexandrov algebra associated to the $C^{\ast}$-algebra $A$. For any morphism $f:A \to B$ of $C^{\ast}$-algebras, we have a morphism $D_{f}:D_{A}\to D_{B}$ of Alexandrov algebras, which sends any element of $D_{A}$ of the form $\{M\in Max(A) \textrm{ | } a\notin M\}$ to the element $\{N\in Max(B) \textrm{ | } f(a)\notin N\}$. This defines a functor D:$\mathbf{C^{\ast}\textrm{-}Alg}\to \textbf{AlexAlg}_{c}$, which is one half of our duality between $C^{\ast}$-algebras and countably compact Alexandrov algebras. 

Conversely, given a countably compact Alexandrov algebra $D$, the set of continuous functions $Max(D)\to {\mathbb R}$ can be bijectively identified with the set ${\cal R}_{D}$ of Alexandrov algebra homomorphisms ${\cal O}({\mathbb R})\to D$ with the property that they send jointly covering families in ${\cal O}({\mathbb R})$ to jointly covering families generated by a countable family of arrows in $D$, that is with the morphisms of sites $({\cal O}({\mathbb R}), J_{{\cal O}({\mathbb R})}^{can})\to (D, C_{D})$. Indeed, since both $Max(D)$ and $\mathbb R$ are sober spaces, the continuous maps $f:Max(D)\to {\mathbb R}$ correspond bijectively with the frame homomorphisms $f^{-1}:{\cal O}({\mathbb R})\to {\cal O}(Max(D))$; but these homomorphisms can be identified with the morphisms of sites $({\cal O}({\mathbb R}), J_{{\cal O}({\mathbb R})}^{can})\to (D, C_{D})$ since they take values in $D\hookrightarrow {\cal O}(Max(D))$ (indeed, $Coz({\cal O}({\mathbb R}))={\cal O}({\mathbb R})$ (cf. \cite{stone}) and $f^{-1}$ sends co-zero sets to co-zero sets). Clearly, any morphism of Alexandrov algebras $g:D\to D'$ yields a morphism of sites $(D, C_{D})\to (D', C_{D'})$ and hence induces a map ${\cal R}_{g}:=g \circ -: {\cal R}_{D}\to {\cal R}_{D'}$. 
Concretely, a morphism of sites $({\cal O}({\mathbb R}), J_{{\cal O}({\mathbb R})}^{can})\to (D, C_{D})$ is a meet-semilattice homomorphism $f:{\cal O}({\mathbb R})\to D$ such that for any family $\{U_{i} \textrm{ | } i\in I\}$ of open sets of $\mathbb R$ there exists a subset $J\subseteq I$ such that $f(J)$ is countable and  $f(\mathbin{\mathop{\textrm{\huge $\vee$}}\limits_{i\in I}}(U_{i}))=\mathbin{\mathop{\textrm{\huge $\vee$}}\limits_{i\in J}}f(U_{i})$.

We can obtain an alternative, more `arithmetic' characterization of the $C^{\ast}$-algebra corresponding to a given countably compact Alexandrov algebra as follows. In section IV 1.1 \cite{stone}, a subcanonical site of definition $(B, C)$ for the topos $\Sh({\mathbb R})$ of sheaves on the topological space $\mathbb R$ is given; we shall show that for any countably compact Alexandrov algebra $D$, the continuous maps $Max(D)\to {\mathbb R}$ can be bijectively identified with the morphisms of sites $(B, C)\to (D, C_{D})$.

The site $(B, C)$ is defined as follows. Let ${\mathbb Q}^{+}$ denote the totally ordered set obtained by adding a top element $\infty$ to the set $\mathbb Q$ of rational numbers, and let ${\mathbb Q}^{-}$ similarly denote ${\mathbb Q}\cup \{\infty\}$. We partially order ${\mathbb Q}^{-}\times {\mathbb Q}^{+}$ by:
\[
(p,q) \leq (p',q') \textrm{ if and only if } p\geq p' \textrm{ and } q\leq q';  
\]  
then ${\mathbb Q}^{-}\times {\mathbb Q}^{+}$ is a meet-semilattice, with top element $(-\infty, \infty)$ and 
$(p, q)\wedge (p', q') = (max\{p, p'\}, min\{q, q'\})$. The Grothendieck topology $C$ on $D$ is generated by the following pullback-stable family of $C$-covering sieves:
\begin{enumerate}[(a)]
\item $\emptyset\in  C(p, q)$ whenever $p\gt q$; 
\item $\{(p, r), (q, s)\} \in C(p, s)$ whenever $p \leq q\lt r \leq s$; and 
\item $\{ (p', q') \textrm{ | } p \lt p' \lt q' \lt q \} \in C(p, q)$ whenever $p\lt q$. 
\end{enumerate}

Notice that $B$ can be identified with the subset of ${\cal O}({\mathbb R})$ consisting of the open sets of the form $(a, \infty)$, $(-\infty, a)$ and $(a, b)$, where $a$ and $b$ are rational numbers, while $C$ can be identified with the Grothendieck topology induced on $B$ by the canonical topology on ${\cal O}({\mathbb R})$.

For any Alexandrov algebra $D$, the morphisms of sites $(B, C)\to (D, C_{D})$ can thus be identified with the meet-semilattice homomorphisms $B\to D$ which send $C$-covering sieves to $C_{D}$-covering sieves, i.e. with the maps $f:{\mathbb Q}^{-}\times {\mathbb Q}^{+} \to D$ such that
\begin{enumerate}[(a)]
\item $f(-\infty, +\infty)=1_{D}$;
\item $f((p, q)\wedge (p',q'))=f(max\{p, p'\}, min\{q, q'\})=f(p,q)\wedge f(p',q')$, for any $(p,q), (p',q')\in {\mathbb Q}^{-}\times {\mathbb Q}^{+}$;
\item $f(p,q)=0_{D}$ whenever $p\gt q$;
\item $f(p,r)\vee f(q,s)=f(p,s)$ whenever $p \leq q\lt r \leq s$;
\item $\mathbin{\mathop{\textrm{\huge $\vee$}}\limits_{p \lt p' \lt q' \lt q}}f(p',q')=f(p,q)$, for any $(p,q) \in {\mathbb Q}^{-}\times {\mathbb Q}^{+}$. 
\end{enumerate}

Let us denote by ${\cal M}_{D}$ the set of such morphisms. As we observed above, the continuous maps $f:Max(D)\to {\mathbb R}$ correspond bijectively with the frame homomorphisms $f^{-1}:{\cal O}({\mathbb R})\to {\cal O}(Max(D))$; but, $C$ being subcanonical, any such homomorphism restricts to a morphism of sites $(B, C)\to (D, C_{D})$ inducing a geometric morphism $\Sh(D, C_{D})\to \Sh(B, C)$ equivalent to $\Sh(f):\Sh(Max(D))\to \Sh({\mathbb R})$. Clearly, any morphism of Alexandrov algebras $g:D\to D'$ yields a morphism of sites $(D, C_{D})\to (D', C_{D'})$ and hence induces a map ${\cal M}_{g}:=g \circ -: {\cal M}_{D}\to {\cal M}_{D'}$. 

Given a countably compact Alexandrov algebra $S$, the $C^{\ast}$-algebra structure on the set of continuous maps $Max(S)\to {\mathbb R}$ can clearly be transferred, via the bijection established above, to a $C^{\ast}$-algebra structure on the set ${\cal M}_{D}$ (resp. on the set ${\cal R}_{D}$), so that the morphisms $D_{f}$ (resp. ${\cal M}_{g}$) become morphisms of the relevant $C^{\ast}$-algebras.
   
Thus we have two functors $D:\mathbf{C^{\ast}\textrm{-}Alg} \to \textbf{AlexAlg}_{c}$ and ${\cal R}:\textbf{AlexAlg}_{c} \to \mathbf{C^{\ast}\textrm{-}Alg}$ (or equivalently, ${\cal M}:\textbf{AlexAlg}_{c} \to \mathbf{C^{\ast}\textrm{-}Alg}$) which are categorical inverses to each other.

Specifically, the following theorem holds.
\begin{theorem}
The functors
\[
D:\mathbf{C^{\ast}\textrm{-}Alg} \to \textbf{AlexAlg}_{c}
\]
and
\[
{\cal R}:\textbf{AlexAlg}_{c} \to \mathbf{C^{\ast}\textrm{-}Alg} 
\]
(or equivalently, 
\[
{\cal M}:\textbf{AlexAlg}_{c} \to \mathbf{C^{\ast}\textrm{-}Alg})
\]
defined above are categorical inverses to each other and yield an equivalence between the category $\mathbf{C^{\ast}\textrm{-}Alg}$ of $C^{\ast}$-algebras and the category $\textbf{AlexAlg}_{c}$ of countably compact Alexandrov algebras.  
\end{theorem}
   
Let us denote by $Max_{r}:\mathbf{C^{\ast}\textrm{-}Alg} \to \textbf{CHaus}$ the maximal spectrum functor, which constitutes one half of Gelfand duality, and by $C:\textbf{CHaus} \to \mathbf{C^{\ast}\textrm{-}Alg}$ the inverse functor associating to a compact Hausdorff space $X$ the $C^{\ast}$-algebra $C(X)$ consisting of all the continuous (bounded) functions $X\to {\mathbb R}$.   

We can represent the dualities which we have established above in the following commutative diagram.

\[  
\xymatrix {
& & \hspace{2.5cm}\textbf{CHaus}^{\textrm{op}} \hspace{2.5cm}   \ar@<0.2ex>[dll]^{C} \ar@<1.6ex>[drr]^{Coz} & &  \\
\mathbf{C^{\ast}\textrm{-}Alg} \ar@<1.6ex>[urr]^{Max_{r}} \ar@<0.6ex>^{D}[rrrr] & & & & \textbf{AlexAlg}_{c} \ar@<0.9ex>[llll]^{{\cal R}\cong {\cal M}} \ar@<0.2ex>[ull]^{Max}}
\]

\subsection{A duality for $T_{1}$ compact spaces}\label{T1}

As another application of the general method of section \ref{generalframework}, we build a different duality between compact Hausdorff spaces and distributive lattices, based on an alternative choice of the normal Wallman bases for the spaces.

First, we notice that if $X$ is a $T_{1}$-space then ${\cal O}(X)$ is a Wallman base for $X$, and for any continuous map $f:X\to Y$ of topological spaces the inverse image $f^{-1}:{\cal O}(Y)\to {\cal O}(X)$ (trivially) restricts to these Wallman bases. Moreover, for any $T_{1}$ compact space $X$, the map $\eta^{X}_{{\cal O}(X)}:X\to Spec({\cal O}(X))$ is surjective on $Max({\cal O}(X))$. Indeed, for any ${\cal M}$ in $Max({\cal O}(X))$, the intersection of all the sets of the form $X\setminus U$ for $U\in {\cal M}$ is non-empty (otherwise, $X$ being compact, $\cal M$ would not be a proper ideal), that is there exists $x\in X$ such that for all $U\in {\cal M}$, $x\notin U$; hence ${\cal M}\subseteq \eta^{X}_{{\cal O}(X)}(x)$ and, by maximality of $\cal M$, $M=\eta^{X}_{{\cal O}(X)}$.

Therefore the category $\textbf{T1Comp}$ of $T_{1}$ compact spaces and continuous maps between them can be identified with a full subcategory of  $\textbf{TopDLat}_{W}$, and Theorem \ref{duality} yields a duality between $\textbf{T1Comp}$ and a full subcategory $\textbf{T1Frm}$ of the category $\textbf{DLat}_{W}$. Specifically, we have a functor
\[
D:\textbf{T1Comp}^{\textrm{op}} \to \textbf{T1Frm}
\]   
sending a topological space $X$ in $\textbf{T1Comp}$ to the open set ${\cal O}(X)$ and a continuous map $f:X\to Y$ of topological spaces in $\textbf{T1Comp}$ to the morphism $f^{-1}:{\cal O}(Y)\to {\cal O}(X)$ in $\textbf{T1Frm}$, and a functor
\[
Max:\textbf{T1Frm} \to \textbf{T1Comp}^{\textrm{op}}
\]   
assigning $Max(D)$ to a distributive lattice $D$ in $\textbf{T1Frm}$ and the continuous map $Max(f):Max(D')\to Max(D)$ to a maximal homomorphism $f:D\to D'$ between lattices in $\textbf{T1Frm}$.

We can describe more explicitly the category $\textbf{T1Frm}$ and the functor $Max$, as follows. The lattices in $\textbf{T1Comp}$ can be characterized as the distributive lattices $D$ such that the canonical morphism $D\to {\cal O}(Max(D))\cong Id_{J_{m}^{D}}(D)$ is an isomorphism. If $D$ is such a lattice then clearly $D$ is a compact frame (i.e., a frame such that every covering of its top elements admits a finite subcovering). Also, it is immediate to see that if $D$ is of the form ${\cal O}(X)$ for a space $X$ such that the map $\eta^{X}_{{\cal O}(X)}:X\to Spec({\cal O}(X))$ is an homeomorphism onto $Max({\cal O}(X))$ then the subtopos 
\[
\Sh(\eta^{X}_{{\cal O}(X)}):\Sh({\cal O}(X), J_{{\cal O}(X)}^{can})=\Sh(X) \to \Sh(Spec({\cal O}(X)))\simeq \Sh({\cal O}(X), J_{{\cal O}(X)}^{coh}),
\]
which is isomorphic to the canonical inclusion 
\[
\Sh({\cal O}(X), J_{{\cal O}(X)}^{can}) \hookrightarrow \Sh({\cal O}(X), J_{{\cal O}(X)}^{coh}),
\]
is isomorphic to the subtopos 
\[
\Sh(Max({\cal O}(X)))\simeq \Sh({\cal O}(X), J_{m}^{{\cal O}(X)})\hookrightarrow \Sh({\cal O}(X), J_{{\cal O}(X)}^{coh}),
\]
from which it follows that $J_{m}^{{\cal O}(X)}=J_{{\cal O}(X)}^{can}$. On the other hand, it is clear that if a distributive lattice $D$ is a frame and $J_{D}^{can}=J_{m}^{D}$ then the canonical morphism $D\to {\cal O}(Max(D))\cong Id_{J_{m}^{D}}(D)$ is an isomorphism. We can thus conclude that the lattices $D$ in $\textbf{T1Frm}$ are precisely the frames $D$ such that $J_{D}^{can}=J_{m}^{D}$; more explicitly, a distributive lattice $D$ satisfies this condition if and only if it is a compact frame (notice that this condition implies that $J_{D}^{can}\subseteq J_{m}^{D}$) which is conjunctive, that is such that for any sieve $S:=\{c_{i}\leq c \textrm{ | } \in I\}$ on $c\in D$ with the property that for every $d\in D$, $c\vee d=1$ implies that there exists a finite subset $J\subseteq I$ such that $\mathbin{\mathop{\textrm{\huge $\vee$}}\limits_{i\in J}}c_{i} \vee d=1$, we have $\mathbin{\mathop{\textrm{\huge $\vee$}}\limits_{i\in J}}c_{i}=c$. 

Summarizing, we have the following result.

\begin{theorem}\label{T1comp}
The functors 
\[
D:\textbf{T1Comp}^{\textrm{op}} \to \textbf{T1Frm}
\]   
and
\[
Max:\textbf{T1Frm} \to \textbf{T1Comp}^{\textrm{op}}
\]   
defined above define a duality between the category $\textbf{T1Comp}$ of $T_{1}$ compact spaces and continuous maps between them and the category $\textbf{T1Frm}$ of compact conjunctive frames and maximal homomorphisms between them.  
\end{theorem}\qed    

We can characterize the maximal ideals of a frame in $\textbf{T1Frm}$ more explicitly, as follows. Given an ideal $I$ of a frame $F$ in $\textbf{T1Frm}$, consider the principal ideal $(\mathbin{\mathop{\textrm{\huge $\vee$}}\limits_{a\in I}}a)$; then we have $I\subseteq (\mathbin{\mathop{\textrm{\huge $\vee$}}\limits_{a\in I}}a)$. If $I$ is proper then, since $F$ is compact, $\mathbin{\mathop{\textrm{\huge $\vee$}}\limits_{a\in I}}a\neq 1$, that is the ideal $(\mathbin{\mathop{\textrm{\huge $\vee$}}\limits_{a\in I}}a)$ is proper; therefore, if $I$ is maximal then $I=(\mathbin{\mathop{\textrm{\huge $\vee$}}\limits_{a\in I}}a)$. The maximal ideals of $F$ can thus be identified with the elements $a$ of $F$ such that the principal ideal $(a)$ is maximal, that is such that for any $a'\in F$, $a\leq a'$ implies $a'=1$; we shall call such elements the \emph{co-atoms} of the frame $F$. Using this identification of the maximal ideals of $F$ with the co-atoms of $F$, we can describe the maximal spectrum $Max(F)$ as the topological space whose underlying set is the set $CoAt(F)$ of co-atoms of $F$ and whose basic open sets are the subsets of the form $\{b\in CoAt(F) \textrm{ | } b\nleq a\}$ for $a\in F$. In these terms the condition for a homomorphism $f:F\to F'$ of frames in $\textbf{T1Frm}$ to be maximal (notice that every morphism in $\textbf{T1Frm}$ is a frame homomorphism) can be expressed as the requirement that for any co-atom $b$ in $F'$ the join $\mathbin{\mathop{\textrm{\huge $\vee$}}\limits_{f(a)\leq b}}a$ should be a co-atom in $F$. 

The following proposition provides an alternative characterizations of conjunctive compact frames.

\begin{proposition}
Let $F$ be a compact frame. Then $F$ is conjunctive if and only if it is co-atomistic in the sense of \cite{maruyama} (i.e., every element of $F$ is the meet of the set of co-atoms greater or equal to it). 
\end{proposition}       

\begin{proofs}
Let us suppose $F$ to be co-atomistic. We want to show that for any sieve $S:=\{c_{i}\leq c \textrm{ | } \in I\}$ on $c\in D$ with the property for every $d\in D$, $c\vee d=1$ implies that there exists a finite subset $J\subseteq I$ such that $\mathbin{\mathop{\textrm{\huge $\vee$}}\limits_{i\in J}}c_{i} \vee d=1$, we have $c=\mathbin{\mathop{\textrm{\huge $\vee$}}\limits_{i\in J}}c_{i}$. Suppose that this latter equality does not hold. Then, $F$ being co-atomistic, there exists a co-atom $a$ of $F$ such that $\mathbin{\mathop{\textrm{\huge $\vee$}}\limits_{i\in J}}c_{i} \leq a$ but $c\nleq a$. Since $a$ is a co-atom we have $a\vee c=1$ and hence there exists a finite subset $J\subseteq I$ such that $\mathbin{\mathop{\textrm{\huge $\vee$}}\limits_{i\in J}}c_{i} \vee a=1$. But since $\mathbin{\mathop{\textrm{\huge $\vee$}}\limits_{i\in J}}c_{i} \leq a$ we have $a=\mathbin{\mathop{\textrm{\huge $\vee$}}\limits_{i\in J}}c_{i} \vee a=1$, which is absurd.

Conversely, suppose that $F$ is conjunctive. We want to prove that for any element $a\in F$, $a$ is equal to the meet in $F$ of all the co-atoms greater or equal to it. Let us denote this meet by $b$; then, clearly, $a\leq b$. To prove that $a=b$ is therefore equivalent to prove that the sieve $S_{a,b}$ generated by the unique arrow $a\leq b$ in $F$ is $J_{F}^{can}$-covering. Since $F$ is conjunctive, $J_{F}^{can}=J_{m}^{F}$ and hence to prove that $S_{a,b}$ is $J_{F}^{can}$-covering amounts precisely to verifying that it is $J_{m}^{F}$-covering, i.e. that for any $c\in F$ such that $c\vee b=1$, $c\vee a=1$. Suppose that this condition does not hold; then there exists $c\in F$ such that $c\vee b=1$ and $c\vee a \neq 1$. Then, by the maximal ideal theorem (and the characterization of maximal ideals as co-atoms in compact frames established above), there exists a co-atom $d$ of $F$ such that $c\vee a \leq d$. Therefore $d$ is a co-atom of $F$ greater or equal to $a$ and hence, by definition of $b$, we have $d\geq b$. Therefore $1=c\vee d=d$, which is absurd.
\end{proofs}

In fact, under this identification between conjunctive compact frames and co-atomistic frames, our duality between the category of $T_{1}$ compact spaces and the category of conjunctive compact frames corresponds precisely to the restriction of the well-known duality between $T_{1}$-spaces and co-atomistic frames (cf. for instance \cite{maruyama}). 
  
Let us now proceed to specialize the duality of Theorem \ref{T1comp} to the context of compact Hausdorff spaces.

Recall that for any $T_{1}$ compact space $X$, $X$ is normal if and only if it is Hausdorff, and $X$ is normal if and only if ${\cal O}(X)$ is a normal lattice (cf. Exercise II 3.6 \cite{stone}). From this remark we immediately deduce that the duality of Theorem \ref{T1comp} restricts to a duality between the category $\textbf{CHaus}$ and the full subcategory $\textbf{T1Frm}_{n}$ of $\textbf{T1Frm}$ on the frames in $\textbf{T1Frm}$ which are normal. We notice that, since every Hausdorff space is sober, every frame homomorphism between frames in $\textbf{T1Frm}$ is maximal. Indeed, by Proposition \ref{maxmorphisms}, the maximal homomorphisms $F\to F'$ between frames $F, F'$ in $\textbf{T1Frm}$ are precisely the morphisms of sites $(F, J_{m}^{F})\to (F', J_{m}^{F'})$; but, $F$ and $F'$ being conjunctive, $J_{m}^{F}=J_{F}^{can}$ and $J_{m}^{F'}=J_{F'}^{can}$, which implies that such morphisms are precisely the frame homomorphisms $F\to F'$. Therefore the category $\textbf{T1Frm}_{n}$ can be described as the full subcategory of the category $\textbf{Frm}$ of frames on the conjunctive (equivalently, co-atomistic) compact normal frames.  

Summarizing, we have the following result.

\begin{theorem}
The functors 
\[
D:\textbf{CHaus}^{\textrm{op}} \to \textbf{T1Frm}_{n}
\]   
and
\[
Max:\textbf{T1Frm}_{n} \to \textbf{CHaus}^{\textrm{op}}
\]   
defined above define a duality between the category $\textbf{CHaus}$ of compact Hausdorff spaces and continuous maps between them and the category $\textbf{T1Frm}_{n}$ of conjunctive (equivalently, co-atomistic) compact normal frames and frame homomorphisms between them. 
\end{theorem}\qed  

Notice that in passing we have established the following fact.

\begin{proposition}
Any frame homomorphism between conjunctive (equivalently, co-atomistic) compact normal frames is maximal.  
\end{proposition}\qed

We mention that in \cite{Ban} a duality different from ours between the category of compact Hausdorff spaces and a subcategory of the category of distributive lattices was established. 

Finally, let us point out an interesting consequence of the duality of Theorem \ref{duality}.

\begin{corollary}
Let $X$ be a $T_{0}$-space and $B$ be a Wallman base of $X$ such that the map $\eta_{B}$ is surjective (equivalently, a homeomorphism) on $Max(B)$ (for example $X$ be compact Hausdorff and $B$ be a normal base for it). Then for any open set $U$ of $X$ belonging to $B$ and any family $\{U_{i} \subseteq U \textrm{ | } i\in I\}$ of open subsets of $U$ belonging to $B$, the union $\mathbin{\mathop{\textrm{\huge $\cup$}}\limits_{i\in I}}U_{i}$ in ${\cal O}(X)$ is equal to $U$ if and only if for any $V$ in $B$ if $V\cup U=X$ there exists a finite subset $J\subseteq I$ such that $\mathbin{\mathop{\textrm{\huge $\cup$}}\limits_{i\in J}}U_{i} \cup V=X$.   
\end{corollary}

\begin{proofs}
The condition in the corollary amounts precisely to the requirement that the maximal topology on $B$ should be equal to the topology on $B$ induced by the canonical topology on ${\cal O}(X)$, and this follows from the duality of Theorem \ref{duality} (specifically, from the definition of the functor $H$). 
\end{proofs}

The following result can be obtained as a particular instance of the corollary above.

\begin{corollary}
Let $X$ be a compact Hausdorff space. Then for any continuous map $f:X\to {\mathbb R}$ and family $\{f_{i}:X\to {\mathbb R} \textrm{ | } i\in I\}$ of continuous maps $f_{i}:X\to {\mathbb R}$ such that for each $i\in I$ $Coz(f_{i})\subseteq Coz(f)$, we have $\mathbin{\mathop{\textrm{\huge $\cup$}}\limits_{i\in I}}Coz(f_{i})=Coz(f)$ if and only if for any continuous map $g:X\to {\mathbb R}$ such that $Coz(f)\cup Coz(g)=X$ there exists a finite subset $J\subseteq I$ such that $\mathbin{\mathop{\textrm{\huge $\cup$}}\limits_{i\in J}}Coz(f_{i}) \cup Coz(g)=X$.    
\end{corollary}\qed

\section{Gelfand spectra}\label{Gelfandspectra}

\subsection{Rings and their reticulations}\label{reticulation}

Given a commutative ring with unit $A$, there is a distributive lattice $L(A)$, called in \cite{HS} the \emph{reticulation} of $A$, such that we have an equivalence
\[
\Sh(L(A), J_{L(A)}^{coh})\simeq \Sh(Spec(A)),
\]
where $Spec(A)$ is the Zariski spectrum of the ring $A$.
 
The lattice $L(A)$ can be described as the coherent syntactic category of the coherent propositional theory ${\mathbb T}^{A}_{p}$ over the signature having a propositional symbol $P_{a}$ for each element $a\in A$, whose axioms are the following:
\[
(\top \vdash P_{1_{A}});
\]
\[
(P_{0_{A}} \vdash \bot);
\]
\[
(P_{a\cdot b} \dashv\vdash P_{a} \wedge P_{b})
\] 
for any $a, b$ in $A$; 
\[
(P_{a + b} \vdash P_{a} \vee P_{b})
\]
for any $a, b \in A$. 

Consider the subspace $Max(A)$ of $Spec(A)$ obtained by inducing on the set of maximal ideals of the ring $A$ the Zariski topology on $Spec(A)$. Recall from \cite{OC11} that the Zariski topology on $Spec(A)$ is homeomorphic, under the complementation map in $\mathscr{P}(A)$, to the subterminal topology on the space of points of the topos $\Sh(L(A), J_{L(A)}^{coh})$.

In \cite{HS} it is argued that the notion of reticulation of a ring can be profitably used for transferring many results about distributive lattices to results about rings and conversely. In this section we shall give a further illustration of this general remark by comparing the maximal spectrum of a commutative ring with unit $A$ with the maximal spectrum of its reticulation.

By the results in \cite{OC11} we have equivalences of toposes
\[
\Sh(L(A), J_{L(A)}^{coh}) \simeq \Sh(Id_{J_{L(A)}^{coh}}(L(A))) \simeq \Sh(X_{L(A)})\simeq \Sh(Rad(A)),
\]
where $X_{L(A)}$ is the space of points of the topos $\Sh(L(A), J_{L(A)}^{coh})$ and $Rad(A)$ is the set of radical ideals of $A$, endowed with the subset-inclusion ordering (cf. Corollary V3.2 \cite{stone}). In fact, there is a further representation of this topos, obtained by cutting the site $(L(A), J_{L(A)}^{coh})$ down to the full subcategory $U$ of $L(A)$ on the objects of the form $[P_{a}]$ (for $a\in A$), where $[P_{a}]$ denotes the equivalence class of the formula $P_{a}$ in $L(A)$; in fact, since finite conjunctions of formulae of the form $P_{a}$ are equivalent in ${\mathbb T}^{A}_{p}$ to formulae of the same form (cf. the third axiom scheme of the theory ${\mathbb T}^{A}_{p}$), any coherent formula over the signature of ${\mathbb T}^{A}_{p}$ is provably equivalent to a finite disjunction of formulae of the form $P_{a}$; hence $U$ is a $J_{L(A)}^{coh}$-dense subcategory of $L(A)$ and the Comparison Lemma yields an equivalence $\Sh(L(A), J_{L(A)}^{coh}) \simeq \Sh(U, J_{L(A)}^{coh}|_{U})$. It is easy to verify that the site $(U, U, J_{L(A)}^{coh}|_{U})$ is categorically equivalent to the site $(S(A), C)$ considered in \cite{OC11}. 

Now, considering the subspace of the closed points (equivalently, of the points which are minimal with respect to the specialization preorder) of the space of points of the topos 
\[
\Sh(X_{L(A)})\simeq \Sh(Spec(A))
\] 
and assuming the prime ideal theorem (for rings), we obtain that this equivalence 
restricts to an equivalence of subtoposes
\[  
\xymatrix {
\Sh(X_{L(A)}) \ar[r]^{\simeq} & \Sh(Spec(A))\\
\Sh(Max(L(A))) \ar[u]  \ar[r]_{\simeq} & \Sh(Max(A)) \ar[u].}
\] 

Assuming the maximal ideal theorem (for distributive lattices), the nucleus $K_{L(A)}$ (cf. section \ref{maxspectrum} for the definition of the nucleus $K_{D}$ for a distributive lattice $D$) on $Id_{J_{L(A)}^{coh}}(L(A))$ corresponds, under the isomorphism 
\[
Id_{J_{L(A)}^{coh}}(L(A)) \cong Rad(A)
\]
to the nucleus $k_{A}$ on $Rad(A)$ defined by the formula
\[
K_{A}(I)=\{a\in A \textrm{ | } (\forall b\in A)(\exists c\in A) (ab+c+abc\in I)\}.
\]  
Indeed, as observed in \cite{almax}, for any (radical) ideal $I$ of $A$ $K_{A}(I)$ coincides with the Jacobson radical of $I$, that is with the intersection of all the maximal ideals of $A$ which contain $I$, and by Proposition 1.6 \cite{almax} for any ideal $V$ of $L(A)$ the ideal $K_{L(A)}(V)$ is the intersection of all the maximal ideals of $L(A)$ which contain $V$; whence our claims follows from the isomorphism of subtoposes established above.

By the duality theorem of \cite{OC6}, the subtopos 
\[
\Sh(Max(A))\hookrightarrow \Sh(Spec(A))
\]
corresponds to a unique quotient of the theory ${\mathbb T}^{A}_{p}$, which we call ${\mathbb T}^{A}_{m}$. In order to give an explicit axiomatization of the theory ${\mathbb T}^{A}_{m}$, we explicitly characterize the Grothendieck topology $M$ on $L(A)$ corresponding to the subtopos 
\[
\Sh(Max(A))\hookrightarrow \Sh(Spec(A))\simeq \Sh(L(A), J_{L(A)}^{coh}).
\]
By the equivalence of subtoposes established above and the equivalence $\Sh(L(A), J_{L(A)}^{coh})\simeq \Sh(U, J_{L(A)}^{coh}|_{U})$, the Grothendieck topology $M$ can be generated by a family of sieves lying in the subcategory $U$.  

Let us apply the characterization of the topology $J_{m}^{D}$ on a distributive lattice $D$ obtained in section \ref{maxspectrum} to the case in which $D$ is the reticulation $L(A)$ of a commutative ring with unit $A$; in this case we denote the topology $J_{m}^{L(A)}$ on $L(A)$ simply by $J_{m}^{A}$. In view of Lemma V 3.2 \cite{stone}, for any sieve $S:=\{[P_{c_{i}}]\leq [P_{c}] \textrm{ | } i\in I\}$ on $[P_{c}]$ in $U$, $S$ is $J_{m}^{A}$-covering if and only if for any finite set of elements $d_{1}, \ldots, d_{n}\in A$ such that $(c, d_{1}, \ldots, d_{n})=A$ there exists a finite subset $J\subseteq I$ such that the ideal generated by the $d_{i}$ (for $i\in \{1, \ldots, n\}$) and the $c_{J_{m}^{B}}$ (for $j\in J$) is the whole of $A$. Notice that for any $a,b\in D$, we have $[P_{a}]\leq [P_{b}]$ in $L(A)$ if and only if there exists an integer $n\gt 0$ and an element $c\in A$ such that $a^{n}=bc$ (cf. for example \cite{OC11}). 

Therefore the theory ${\mathbb T}^{A}_{m}$ can be axiomatized by adding to the axioms of ${\mathbb T}^{A}_{p}$ all the sequents of the form 
\[
(P_{c} \vdash \mathbin{\mathop{\textrm{\huge $\vee$}}\limits_{i\in I}}P_{c_{i}})
\]
for any elements $\{c_{i} \textrm{ | } i\in I\}$ and $c$ of $A$ such that for any $i\in I$ there exists an integer $n_{i}\gt 0$ and an element $u_{i}\in A$ such that $c_{i}^{n_{i}}=c u_{i}$ and for any finite set of elements $d_{1}, \ldots, d_{n}\in A$ such that $(c, d_{1}, \ldots, d_{n})=A$ there exists a finite subset $J\subseteq I$ such that the ideal generated by the $d_{i}$ (for $i\in \{1, \ldots, n\}$) and the $c_{J_{m}^{B}}$ (for $j\in J$) is the whole of $A$.

The following result represents the analogue for rings of Theorem \ref{dlatmax} above. 

\begin{theorem}\label{sobermax}
Let $A$ be a commutative ring with unit such that the space $Max(A)$ is sober (for example, $A$ is a $C^{\ast}$-algebra). Then a prime ideal $P$ of $A$ is maximal if and only if for any elements $\{c_{i} \textrm{ | } i\in I\}$ and $c$ of $A$ with the property that 
\begin{enumerate}
\item for any $i\in I$ there exists an integer $n_{i}\gt 0$ and an element $u_{i}\in A$ such that $c_{i}^{n_{i}}=c u_{i}$ and 
\item for any finite set of elements $d_{1}, \ldots, d_{n}\in A$ such that $(c, d_{1}, \ldots, d_{n})=A$ there exists a finite subset $J\subseteq I$ such that the ideal generated by the $d_{i}$ (for $i\in \{1, \ldots, n\}$) and the $c_{j}$ (for $j\in J$) is the whole of $A$, 
\end{enumerate}
if $c\notin P$ then $c_{i}\notin P$ for some $i\in I$.  
\end{theorem}\qed

We shall call a prime ideal $P$ of a commutative ring with unit $A$ satisfying the condition in the statement of Theorem \ref{sobermax} an \emph{almost maximal} ideal of $A$.

Let us now proceed to establish a corresponding characterization, holding under the assumption of the maximal ideal theorem, of the homomorphisms of commutative rings with unit $f:A\to B$ such that $f^{-1}:Spec(B)\to Spec(A)$ restricts to a (continuous) map $Max(B)\to Max(A)$. 

By applying Proposition \ref{maxmorphisms} to the case of the distributive lattice homomorphism $L(f):L(A)\to L(B)$ induced by the homomorphism $f:A\to B$, we obtain the following result, which represents the ring-theoretic analogue of it.

\begin{theorem}\label{morphisms}
Let $f:A\to B$ be a homomorphism of commutative rings with unit $A$ and $B$. If the continuous map $f^{-1}:Spec(B)\to Spec(A)$ restricts to a (continuous) map $Max(B)\to Max(A)$ then for any elements $\{c_{i} \textrm{ | } i\in I\}$ and $c$ of $A$ with the property that for any $i\in I$ there exists an integer $n_{i}\gt 0$ and an element $u_{i}\in A$ such that $c_{i}^{n_{i}}=c u_{i}$ and for any finite set of elements $d_{1}, \ldots, d_{n}\in A$ such that $(c, d_{1}, \ldots, d_{n})=A$ there exists a finite subset $J\subseteq I$ such that the ideal generated by the $d_{i}$ (for $i\in \{1, \ldots, n\}$) and the $c_{j}$ (for $j\in J$) is the whole of $A$, for any finite set of elements $d_{1}', \ldots, d_{m}'\in B$ such that $(f(c), d_{1}', \ldots, d_{m}')=B$ there exists a finite subset $K\subseteq I$ such that the ideal generated by the $d_{i}'$ (for $i\in \{1, \ldots, m\}$) and the $f(c_{i})$ (for $i\in K$) is the whole of $B$. The converse implication holds if $Max(A)$ is sober. 
\end{theorem}\qed

We shall say that a homomorphism $f:A\to B$ of commutative rings with unit is \emph{maximal} if $f^{-1}:Spec(B)\to Spec(A)$ restricts to a (continuous) map $Max(B)\to Max(A)$.

\subsection{The case of commutative $C^{\ast}$-algebras}

Recall that there are two kinds of Gelfand duality: a real version and a complex one. The real version, due to Stone \cite{Stone} and described in \cite{stone}, gives a duality between a full subcategory of the category of commutative rings with unit, which we call the category of \emph{real} $C^{\ast}$-algebras, and the category of compact Hausdorff spaces, while the complex version, also known as the classical Gelfand duality, gives a duality between the category of (complex) commutative $C^{\ast}$-algebras and the category of compact Hausdorff spaces. 

Recall that if $A$ and $B$ are real $C^{\ast}$-algebras then a $C^{\ast}$-algebra homomorphism $A\to B$ is defined simply as a homomorphism $A\to B$ of commutative rings with unit, while if $A$ and $B$ are complex $C^{\ast}$-algebras a $C^{\ast}$-algebra homomorphism $A\to B$ is defined as a ring homomorphism $A\to B$ which commutes with the involution $\ast$.

Before proceeding further let us remark some useful facts about Gelfand duality. In the real case, one half of the duality sends a real $C^{\ast}$-algebra $A$ to the maximal spectrum $Max(A)$ of $A$ (considered as a ring). The maximal spectrum $Max(A)$ can be alternatively be described as the set $C(A)$ of non-zero ${\mathbb R}$-algebra homomorphisms $A\to {\mathbb R}$, endowed with the weak $*$-topology. In fact, for any non-zero ${\mathbb R}$-algebra homomorphisms $f:A\to {\mathbb R}$, $f^{-1}(0)$ is a maximal ideal of $A$ (since for any $b$ such that $b\notin f^{-1}(0)$ there exists a real number $r$ such that $rf(b)=1$ and hence $(r\cdot 1_{A})b -1\in f^{-1}(0)$), while for any maximal ideal $M$ of $A$ the quotient $A\slash M$ is isomorphic to $\mathbb R$ (cf. Theorem IV 4.7 \cite{stone}). In the complex case there is an analogous characterization of the maximal ideals as the $\mathbb C$-algebra homomorphisms from the algebra to the field $\mathbb C$ of complex numbers. Actually, these two different descriptions of the maximal spectrum of a $C^{\ast}$-algebra can be interpreted as a Morita-equivalence between two propositional theories over different signatures, one of `algebraic' nature whose models in $\Set$ are the maximal ideals on the $C^{\ast}$-algebra and one of `analytic' nature, whose models in $\Set$ are the algebra homomorphism from the $C^{\ast}$-algebra to $\mathbb R$ or $\mathbb C$ (cf. section \ref{toposintgelfand} below).      

We shall now apply the characterizations of the maximal ideals of ring whose maximal spectrum is sober given by Theorem \ref{sobermax} and the characterizations of maximal homomorphisms of rings given by Theorem \ref{morphisms} in the context of $C^{\ast}$-algebras and Gelfand duality. 

\begin{theorem}
\begin{enumerate}[(i)]
\item Let $A$ be a (real or complex) $C^{\ast}$-algebra. Then for any prime ideal $P$ of $A$, $P$ is maximal if and only if it is almost maximal;

\item Let $f:A\to B$ a homomorphisms of rings between (real or complex) $C^{\ast}$-algebras $A$ and $B$. Then $f$ is a maximal homomorphism of rings $A\to B$.
\end{enumerate}
\end{theorem}

\begin{proofs}
$(i)$ In order to apply Theorem \ref{sobermax}, it suffices to show that if $A$ is a $C^{\ast}$-algebra (whether real or complex) then the space $Max(A)$ is sober, and this immediately follows from Gelfand duality since every Hausdorff space is sober.

$(ii)$ We shall prove the result for real $C^{\ast}$-algebras, the complex case being entirely analogous to it. In order to apply Theorem \ref{morphisms}, we have to verify that for any homomorphism of real $C^{\ast}$-algebras $f:A\to B$ satisfies the property that the inverse image map $f^{-1}:Spec(B)\to Spec(A)$ restricts to a function $Max(B)\to Max(A)$. We recall that one half of Gelfand duality sends a homomorphism of $C^{\ast}$-algebras $f:A\to B$ to the continuous map $Max(B)\to Max(A)$ which corresponds, under the identifications $Max(B)\cong C(B)$ and $Max(A)\cong C(A)$, to the map $C(B)\to C(A)$ sending a $\mathbb R$-algebra homomorphism $g:B\to {\mathbb R}$ to the composite $g\circ f:A\to {\mathbb R}$. But under the identifications $Max(B)\cong C(B)$ and $Max(A)\cong C(A)$ mentioned above, this assignment is immediately seen to correspond exactly to taking the inverse image of maximal ideals under $f$, as required. We can thus appeal to Theorem \ref{morphisms} to conclude our thesis. 
\end{proofs}

Let us define a commutative ring with unit $A$ to be \emph{conjunctive} if its reticulation $L(A)$ is conjunctive as a distributive lattice (in the sense of section \ref{maxspectrum} above). In order to get an explicit characterization of conjunctive rings, we shall need the following lemma about subcanonical topologies.

\begin{lemma}
Let $({\cal C}, J)$ be a subcanonical Grothendieck site and $\cal D$ be a full $J$-dense subcategory of $C$. Then for any Grothendieck topology $K$ on $\cal C$ which contains $J$, $K$ is subcanonical if and only if any $K|_{D}$-covering sieve on an object $d\in {\cal D}$ generates an effective-epimorphic sieve in $\cal C$ (i.e., the sieve $R$ in $\cal C$ generated by it forms a colimit cone under the (possibly large) diagram consisting of the domains of all the morphisms in $R$, and all the morphisms over $d$ between them).   
\end{lemma}

\begin{proofs}
The `only if' implication is clear so it remains to prove the `if' one. 

Let $R$ be a $K$-covering sieve on an object $c\in {\cal C}$. We want to show that $R$ is effective epimorphic, that is for any object $c'\in {\cal C}$ and cone $\{s_{f}:dom(f) \to c' \textrm{ | } f\in R\}$ with vertex $c'$ over the diagram in $\cal C$ formed by the objects of the form $dom(f)$ (for $f\in R$) and the arrows over $c'$ between them there exists a unique arrow $\xi:c\to c'$ in $\cal C$ such that $\xi\circ f=s_{f}$ for every $f\in R$.

Since $\cal D$ is $J$-dense in $\cal C$, there exists a $J$-covering sieve $Z$ on $c$ generated by arrows whose domains are in $\cal D$. Notice that, since by our hypothesis $J$ is subcanonical, the sieve $Z$ is effective-epimorphic (in the sense of the definition at p. 542 \cite{El}). Consider the sieves $f^{\ast}(R)\cap {\cal D}\in K|_{\cal D}(dom(f))$ (for $f\in Z$). For any $f\in Z$, the arrows $\{s_{f\circ h} \textrm{ | } h\in f^{\ast}(R)\cap {\cal D}\}$ form a cone with vertex $c'$ over the diagram formed by the objects of the form $dom(h)$ (for $h\in f^{\ast}(R)\cap {\cal D}$) and the arrows over $dom(f)$ between them. By our hypotheses we can thus conclude that there exists a unique arrow $t_{f}:dom(f)\to c'$ such that $t_{f}\circ h=s_{f\circ h}$ for every $h\in f^{\ast}(R)\cap {\cal D}$. Now, consider the set of arrows $\{t_{f}:dom(f)\to c' \textrm{ | } f\in Z\}$. Let us show that they form a cone with vertex $c'$ on the objects of the form $dom(f)$ (for $f\in Z$) and the arrows in $\cal C$ over $c$ between them; we have to verify that for any arrows $f,f'\in Z$ and arrow $u:dom(f)\to dom(f')$ in $\cal C$ such that $f'\circ u=f$, $t_{f'}\circ u=t_{f}$. But this equality follows from the fact that the sieve $f^{\ast}(R)\cap {\cal D}$ is effective-epimorphic in $\cal C$ since for every $h\in f^{\ast}(R)\cap {\cal D}$, $t_{f'}\circ u\circ h=s_{f\circ h}$; indeed, $t_{f'}\circ u\circ h=t_{f'}\circ (u\circ h)=s_{f'\circ u\circ h}=s_{f\circ h}$. Now, $Z$ being effective-epimorphic, there is a unique arrow $\xi:c\to c'$ such that $\xi\circ f=t_{f}$ for any $f\in Z$. To conclude our proof, it will be enough to verify that for any $g\in R$, $\xi \circ g=s_{g}$. Since $Z$ is universally effective-epimorphic (since it belongs to the subcanonical topology $J$), the sieve $g^{\ast}(Z)$ is effective-epimorphic and hence $\xi \circ g=s_{g}$ if and only if for any $u\in g^{\ast}(Z)$, $\xi \circ g\circ u=s_{g}\circ u$. But $g\circ u\in Z$ whence $\xi\circ g\circ u = t_{g\circ u}$, and we have $t_{g\circ u}=s_{g}\circ u$ since for any $h\in (g\circ u)^{\ast}(R)\cap {\cal D}$ $t_{g\circ u}\circ h=s_{g\circ u\circ h}=s_{g} \circ u\circ h$ (the latter equality holding since $g\in R$ and the family $\{s_{f}:dom(f) \to c' \textrm{ | } f\in R\}$ is a cone over $R$).      
\end{proofs}

We can apply the lemma in the context of the site $(L(A), J_{L(A)}^{coh})$, with $J$ being equal to $J_{L(A)}^{coh}$, $K$ being equal to the maximal topology $J^{L(A)}_{m}$ and $\cal D$ being equal to the $J_{L(A)}^{coh}$-dense subcategory $U$ of $L(A)$. Recalling Lemma V 3.2 \cite{stone} and the fact that every object of $L(A)$ can be expressed as a finite join of elements in $U$ we easily obtain the following characterization result.

\begin{proposition}
A commutative ring with unit $A$ is conjunctive if and only if for any elements $\{c_{i} \textrm{ | } i\in I\}$ and $c$ of $A$ with the property that for any $i\in I$ there exists an integer $n_{i}\gt 0$ and an element $u_{i}\in A$ such that $c_{i}^{n_{i}}=c u_{i}$ and for any finite set of elements $d_{1}, \ldots, d_{n}\in A$ such that $(c, d_{1}, \ldots, d_{n})=A$ there exists a finite subset $J\subseteq I$ such that the ideal generated by the $d_{i}$ (for $i\in \{1, \ldots, n\}$) and the $c_{j}$ (for $j\in J$) is the whole of $A$, for any elements $a_{1}, \ldots, a_{k}$ of $A$ such that for every $i\in I$ a power of $c_{i}$ belongs to the ideal $(a_{1}, \ldots, a_{k})$ generated by the $a_{i}$, a power of $c$ belongs to $(a_{1}, \ldots, a_{k})$. 
\end{proposition}\qed

Notice that this proposition can be profitably applied in connection with the result of M. Hochster (cf. \cite{hochster}) that every distributive lattice is, up to isomorphism, of the form $L(A)$ for some ring $A$.

If $Max(A)$ is sober then the condition for $A$ to be conjunctive radically simplifies. Specifically, we have
the following criterion.

\begin{proposition}
Let $A$ be a commutative ring with unit such that the space $Max(A)$ is sober (for example, a commutative $C^{\ast}$-algebra). Then $A$ is conjunctive if and only if for any elements $a,b\in A$, if $a$ and $b$ are contained in the same maximal ideals then $a$ and $b$ are contained in the same prime ideals.    
\end{proposition}

\begin{proofs}
It suffices to observe that if $Max(A)$ is sober then the maximal ideals of $A$ can be identified with the points of the topos $\Sh(Max(A))\simeq \Sh(L(A), J_{m}^{L(A)})$. But $J_{m}^{L(A)}$ is subcanonical if and only if for any $a,b\in A$, $[P_{a}]\neq [P_{b}]$ (equivalently, there is a prime ideal $P$ such that $a\in P$ and $b\notin P$ or $b\in P$ and $a\notin P$) implies $a_{m}([P_{a}])\neq a_{m}([P_{b}])$ where $a_{m}:\Sh(L(A), J_{L(A)}^{coh})\to \Sh(L(A), J_{m}^{L(A)})$ is the associated sheaf functor (equivalently, since the maximal ideals of $A$ are a separating of points of the topos $\Sh(L(A), J_{m}^{L(A)})$, there is a maximal ideal $M$ of $A$ such that $a\in M$ and $b\notin M$ or $b\in M$ and $a\notin M$).
\end{proofs}

\begin{remark}
\begin{enumerate}[(a)]
\item If $A$ is semisimple, i.e. its Jacoson radical is the zero ideal and $Max(A)$ is sober (for example if $A$ is a finite dimensional $C^{\ast}$-algebra) then $A$ is conjunctive. Indeed, the Jacobson radical of a ring can be characterized as the intersection of all the maximal ideals of the ring, and for any elements $a, b\in A$, $a$ and $b$ are contained in the same prime (resp. maximal) ideals of $A$ if and only if $a-b$ is contained in every prime (resp. maximal) ideal of $A$. 

\item If $A$ has the property that every prime ideal is equal to its Jacobson radical (i.e., to the intersection of all the maximal ideals containing it) then $A$ is conjunctive; for example, any Boolean ring is conjunctive. Recall from \cite{almax} that the Jacobson radical of an ideal $I$ can be described as the set of all the elements $a$ of $A$ such that $1-ab$ is invertible modulo $I$ for every $b\in B$; from this it immediately follows that if a ring $A$ satisfies the property that for any prime ideal $P$ of $A$ and any element $a\in A$, if $1-ab$ is invertible modulo $P$ for every $b\in A$ then $a\in P$ then $A$ is conjunctive.
\end{enumerate}
\end{remark}

\begin{proposition}
A commutative $C^{\ast}$-algebra $A$ is conjunctive if and only if the distributive lattice $L(A)$ is isomorphic to the lattice $Coz(Max(A))$ under the canonical homomorphism $L(A)\to Coz(Max(A))$. 
\end{proposition}

\begin{proofs}
If $A$ is a $C^{\ast}$-algebra then by Gelfand duality the factorization $\xi_{A}:L(A)\to {\cal O}(Max(A))$ of the map $A\to {\cal O}(Max(A))$ sending an element $a\in A$ to the open set $\{M\in Max(A) \textrm{ | } a\notin M\}$ through the canonical map $A\to L(A)$ takes values in $Coz(Max(A))$ and is equal to the composite of the inclusion $Coz(Max(A))\hookrightarrow {\cal O}(Max(A))$ with the canonical homomorphism $L(A)\to Coz(Max(A))$. Therefore, this latter morphism is an isomorphism if and only if the map $\xi_{A}$ is injective; but $L(A)$ is conjunctive if and only if the map $L(A)\to {\cal O}(Max(L(A)))$ sending an element $d\in L(A)$ to the open set $\{M\in Max(L(A)) \textrm{ | } d\in M\}$ is injective, and this map corresponds exactly to the map $\xi_{A}$ under the isomorphism ${\cal O}(Max(L(A)))\cong {\cal O}(Max(A))$. 
\end{proofs}

\subsection{The topos-theoretic interpretation}\label{toposintgelfand}

In this section we point out several Morita-equivalences which naturally arise in the context of Gelfand duality for $C^{\ast}$-algebras. These equivalences are important in that they formalize the different approaches to the construction of Gelfand spectra, and allow an effective transfer of information between them according to the methodologies introduced in \cite{OC10}; specifically, each way of constructing the spectrum of a $C^{\ast}$-algebra corresponds to a different site of definition for the topos of sheaves on it. Actually, this situation is analogous to that of the different constructions of the Zariski spectrum of a ring, which we interpreted in \cite{OC11} as a collection of Morita-equivalences, as well as to the different ways for building the the Stone-\v{C}ech compactification of a locale, which we interpreted in section \ref{WSC} above. In fact, it is often the case that different points of view on a given mathematical object or different ways for constructing it, as well as mathematical equivalences or dualities of various nature can be formalized as Morita-equivalences (cf. \cite{OC10} for a comprehensive discussion of these issues). 

The construction of the spectrum $X$ of a given $C^{\ast}$-algebra $A$ in terms of the lattice $Coz(X)$ of co-zero sets on it can be formalized by the Morita-equivalence 
\[
\Sh(X)\simeq \Sh(Coz(X), J_{{\cal O}(X)}^{can}|_{Coz(X)}) \simeq \Sh(D_{A}, C_{D_{A}}),
\]
where $D_{A}$ is the Alexandrov algebra associated to the $C^{\ast}$-algebra $A$ as in section \ref{dualityalexandrovalg} and $C_{D_{A}}$ is the countable topology on it.

The representation of $X$ as the space $Max(A)$ of maximal ideals of the $C^{\ast}$-algebra $A$ in the Zariski topology can be expressed by the Morita-equivalence
\[
\Sh(X)\simeq \Set[{\mathbb T}^{A}_{m}],
\]
where $\Set[{\mathbb T}^{A}_{m}]$ is the classifying topos for the theory ${\mathbb T}^{A}_{m}$ of maximal ideals on $A$ defined in section \ref{reticulation} above.

The representation of $A$ as the space of multiplicative linear functionals on $A$ in the weak$^{\ast}$ topology can be expressed as a Morita-equivalence
\[
\Sh(X)\simeq \Set[{\mathbb M}Fn A],
\]
where $\Set[{\mathbb M}Fn A]$ is the classifying topos of the propositional theory ${\mathbb M}Fn A$ of multiplicative linear functionals on $A$ introduced in \cite{BM1}.   

The Gelfand-Mazur isomorphism can thus be interpreted as a Morita-equivalence
\[
\Set[{\mathbb T}^{A}_{m}]\simeq  \Set[{\mathbb M}Fn A].
\]

Also, we have representations of $\Sh(X)$ as subtoposes of $\Sh(L(A), J_{L(A)}^{coh})$ and of $\Sh(U, J_{L(A)}^{coh}|U)$ (cf. section \ref{reticulation} above for the definition of the category $U$):
\[
\Sh(X)\simeq \Sh(L(A), J_{m}^{L(A)})\simeq \Sh(U, J_{m}^{L(A)}|_{U}). 
\]

These Morita-equivalences capture one half of Gelfand duality, namely the representation of the spectrum of a $C^{\ast}$-algebra in terms of the algebra itself; note that the other half of the duality, namely the construction of the $C^{\ast}$-algebra $A$ corresponding to a given compact Hausdorff space $X$ also admits a natural topos-theoretic interpretation, as a sheaf representation result of $A$ as the ring of global sections of a sheaf of local rings defined on $X$ (cf. Corollary V 3.8 \cite{stone}).

\vspace{0.4cm}

{\bf Acknowledgements.} I am grateful to Vincenzo Marra for bringing my attention to \cite{Ban}.

\end{document}